\documentclass[journal,twoside,web]{ieeecolor}
\usepackage{amsmath,amssymb,amsfonts}

\def\logoname{LOGO-generic-web}
\definecolor{subsectioncolor}{rgb}{0,0.541,0.855}
\def\journalname{IEEE Transactions on Automatic Control}

\def\BibTeX{{\rm B\kern-.05em{\sc i\kern-.025em b}\kern-.08em
    T\kern-.1667em\lower.7ex\hbox{E}\kern-.125emX}}

\usepackage[T1]{fontenc}
\usepackage{comment}

\usepackage{tikz}
\usepackage{pgfplots}
\usepackage{subfig}
\usepackage{caption}
\usepackage{siunitx}

\newcommand{\cN}{\mathcal{N}}

\newcommand{\cG}{\mathcal{G}}

\newcommand{\cD}{\mathcal{D}}

\newcommand{\cF}{\mathcal{F}}

\usetikzlibrary{arrows}
\usetikzlibrary{fit}\usetikzlibrary{calc}
  \pgfdeclarelayer{background}
  \pgfsetlayers{background,main}

\newtheorem{Prop}{Proposition}[section]

\newtheorem{Thm}[Prop]{Theorem}
\newtheorem{Cor}[Prop]{Corollary}

\newtheorem{Def}[Prop]{Definition}
\newtheorem{Rem}[Prop]{Remark}

\newcommand{\RN}[1]{\uppercase\expandafter{\romannumeral#1}}
\newcommand{\eps}{\varepsilon}

\newcommand{\N}{{\mathbb{N}}}

\newcommand{\R}{{\mathbb{R}}}

\newcommand{\cTT}{{\mathbb T}^{n,q}_h}
\newcommand{\cCC}{C([-h,\infty),\R^n)}
\newcommand{\cLL}{L_{\rm loc}^\infty(\R_{\ge 0},\R^q)}

\DeclareMathOperator{\esssup}{ess\, sup}

\newcommand{\setdef}[2]{\left\{\, #1 \left|\, \vphantom{#1} #2\right.\right\}}
\newcommand{\ddt}{\tfrac{\text{\normalfont d}}{\text{\normalfont d}t}}

\DeclareMathOperator{\sgn}{sgn}
\DeclareMathOperator{\loc}{loc}

\DeclareMathOperator{\sat}{sat}

{\left(\begin{smallmatrix}}
{\end{smallmatrix}\right)}

{\left[\begin{smallmatrix}}
{\end{smallmatrix}\right]}

\renewcommand{\theenumi}{{\rm(\roman{enumi})}} 

\begin{document}

\title{Input-constrained funnel control\\ of nonlinear systems}
\author{Thomas Berger\thanks{Funded by the Deutsche Forschungsgemeinschaft (DFG, German
Research Foundation) -- Project-ID 471539468.}
\thanks{Thomas Berger is with the Universit\"at Paderborn, Institut f\"ur Mathematik, Warburger Str.~100, 33098~Paderborn, Germany (e-mail: thomas.berger@math.upb.de). }}

\maketitle

\begin{abstract}
We study tracking control for uncertain nonlinear multi-input, multi-output systems modelled by $r$-th order functional differential equations (encompassing systems with arbitrary strict relative degree) in the presence of input constraints. The objective is to guarantee the evolution of the tracking error within a performance funnel with prescribed asymptotic shape (thus achieving desired transient and asymptotic accuracy objectives), for any sufficiently smooth reference signal. We design a novel funnel controller which, in order to satisfy the input constraints, contains a dynamic component which widens the funnel boundary whenever the input saturation is active. This design is model-free, of low-complexity and extends earlier funnel control approaches. We present a simulation where the controller is compared to these approaches.
\end{abstract}

\begin{IEEEkeywords}
adaptive control,
functional differential equations,
funnel control,
input constraints,
nonlinear systems.
\end{IEEEkeywords}


%
\section{Introduction}\label{Sec:Intr}
%

\IEEEPARstart{W}{e} study funnel control for the class of nonlinear systems modelled by the $r$-th order functional differential equation
\begin{equation}\label{eq:Sys}
    \begin{aligned}
      y^{(r)}(t) &= f\big(d(t), T(y,\dot y,\ldots,y^{(r-1)})(t),u(t)\big),\\
      y|_{[-h,0]} &= y^0 \in C^{r-1}([-h,0],\R^m),
    \end{aligned}
\end{equation}
with unknown nonlinear function $f\in C(\R^p\times \R^q\times\R^m,\R^m)$ and unknown operator~$T$ which satisfy a sector bound property (see Section~\ref{Ssec:SysClass}), unknown bounded disturbance~$d$ and unknown initial trajectory~$y^0$ in the presence of input constraints
\begin{equation}\label{eq:IC}
u(t) = \sat(v(t))
\end{equation}
with \textit{known} saturation function $\sat$ and control function~$v$ provided by the to-be-designed controller. Since the objective of funnel control is to achieve a prescribed performance of the tracking error, that is $\|y(t) - y_{\rm ref}(t)\| < \psi(t)$ for some given reference signal $y_{\rm ref}$ and funnel function~$\psi$, a conflict of objectives arises: It is not possible to simultaneously satisfy the input and output constraints for any given bounded reference signal.

In this paper, we consider the input constraints to be \textit{hard constraints}, being imposed by the physical limitations of the system. On the other hand, the output constraints are considered to be \textit{soft constraints}, which can be weakened whenever this is inevitable in order to meet the input constraints. To achieve this, we propose a novel control design, where the funnel function~$\psi(t)$ is no longer prescribed for all $t\ge 0$ as in classical funnel control (see e.g.~\cite{BergIlch21,IlchRyan02b}), but it is dynamically generated and becomes part of the controller design. The generation mechanism for $\psi(t)$ is such that it has a prescribed shape (determined by the parameters in the differential equation which can be chosen a priori by the designer) whenever the saturation is not active, that is $u(t) = v(t)$. In this case, the controller satisfies the input constraints imposed by~\eqref{eq:IC} and achieves the prescribed performance of the tracking error; it further exhibits the same controller performance as the funnel controllers proposed in~\cite{BergIlch21,BergLe18a}. When the saturation is active the performance funnel described by~$\psi(t)$ is widened according to a dynamic equation so that the input constraints are still met~-- in this case, it deviates from the prescribed shape. As soon as the saturation becomes inactive again, the performance funnel recovers its desired shape exponentially fast.

The concept of funnel control was developed in the seminal work~\cite{IlchRyan02b} (see also the recent survey in~\cite{BergIlch21}) and proved advantageous in a variety of applications such as control of industrial servo-systems~\cite{Hack17} and underactuated multibody systems~\cite{BergDrue21}, control of electrical circuits~\cite{SenfPaug14}, control of peak inspiratory pressure~\cite{PompWeye15}, adaptive cruise control~\cite{BergRaue20} and even the control of infinite-dimensional systems such as a boundary controlled heat equation~\cite{ReisSeli15b} and a moving water tank~\cite{BergPuch22}. 

Funnel control with input saturation was first investigated in~\cite{IlchTren04} for the specific application of chemical reactor models and in a more general context in~\cite{HopfIlch10a,HopfIlch10b} for systems with relative degree one and in~\cite{HackHopf13} for systems with relative degree two; this approach has been applied to funnel control with anti-windup for synchronous machines in~\cite{Hack15a}. However, in the aforementioned works it was simply shown that classical funnel control is feasible for a sufficiently large saturation level~-- in the present paper this level can be arbitrarily small. Another approach to funnel control with guaranteed input constraints is bang-bang funnel control, where the control signal switches between only two values. This approach was introduced in~\cite{LibeTren13b} for (undisturbed) nonlinear single-input, single-output systems with arbitrary relative degree. However, the bang-bang funnel control design requires various complicated feasibility assumptions and in particular the two control values must be sufficiently large (typically much larger than actually needed).

A relative of funnel control is prescribed performance control, developed in~\cite{BechRovi08}, see also the important work~\cite{BechRovi14} where the complexity issue of this approach has been solved. The problem of input constraints has been addressed within this approach e.g.\ in~\cite{LiXian18}, where neural networks are used to approximate the nonlinearities, and in~\cite{ChenZhan19}, where additionally a neural observer is incorporated in the controller design. However, such approximation techniques can drastically increase the controller complexity and are avoided in the present paper. In the work~\cite{WangHu19} no approximations are needed (and hence the controller is of low complexity), however the proof contains an error and simulations also show that the proposed controller is infeasible in general. The problem is that the scaling parameter~$\kappa$ in the $\chi$-dynamics is chosen as a constant, but actually it needs to depend on the input.


Funnel control for systems with arbitrary relative degree was considered in~\cite{BergIlch21,BergLe18a}. The novel input-constrained funnel control design that we propose in this paper extends these approaches in the following aspects:

\begin{itemize}
  \item Compared to~\cite{BergLe18a} a much more general class of systems is allowed here, similar to~\cite{BergIlch21}. However, we do not require the restrictive high-gain property of the nonlinearity~$f$ or the minimum phase property (characterized by a BIBO property of the operator~$T$) imposed in~\cite{BergIlch21}. On the other hand, we require a sector bound property of~$f$ and~$T$. This condition cannot be dispensed in general, because of the input saturation, see Remark~\ref{Rem:sector} for more details.
  \item The new controller is able to handle arbitrary input constraints~\eqref{eq:IC}. Even if the saturation is never active, i.e., $u(t) = v(t)$ for all $t\ge 0$ for any solution of the closed-loop system, then the new controller is able to guarantee a prescribed performance of the tracking error as in~\cite{BergIlch21,BergLe18a}, with exponentially decaying funnel boundaries.
\end{itemize}

\subsection{Nomenclature}\label{Ssec:Nomencl}

In the following let $\N$ denote the natural numbers, $\N_0 = \N \cup\{0\}$, and $\R_{\ge 0} =[0,\infty)$. By $\|x\|$ we denote the Euclidean norm of $x\in\R^n$. For some interval $I\subseteq\R$, some $V\subseteq\R^m$ and $k\in\N$, $L^\infty(I, \R^{n})$ $\big(L^\infty_{\rm loc} (I, \R^{n})\big)$ is the Lebesgue space of measurable, (locally) essentially bounded {functions} $f\colon I\to\R^n$, $W^{k,\infty}(I,  \R^{n})$ is the Sobolev space of all functions
$f:I\to\R^n$ with $k$-th order weak derivative $f^{(k)}$ and $f,f^{(1)},\ldots,f^{(k)}\in L^\infty(I, \R^{n})$, and
 $C^k(V,  \R^{n})$ is the set of  $k$-times continuously differentiable functions  $f:  V  \to \R^{n}$, with $C(V,  \R^{n}) := C^0(V,  \R^{n})$.

\subsection{System Class}\label{Ssec:SysClass}

We consider functional differential equations of the form~\eqref{eq:Sys} incorporating an operator~$T$ of the following class.

\begin{Def}
For $n,q\in\N$ and $h \geq 0$ the set $\cTT$ denotes the class of operators $T\colon \cCC\to \cLL$ with the following properties.
\begin{enumerate}
\item[\textbf{(P1)}]
 $T$ is causal, i.e., for all $\zeta$, $\xi \in \cCC$ and all $t\ge 0$,
\[
\zeta|_{[-h,t]} =\xi|_{[-h,t]} ~~\implies~~ T(\zeta)|_{[0,t]}=T(\xi)|_{[0,t]}.
\]
\item[\textbf{(P2)}]
 $T$ is locally Lipschitz, i.e., for each $t\ge 0$ and all $\xi\in C([-h,t],\R^{n})$, there exist positive constants $c_0, \delta, \tau >0$ such that, for all $\zeta_1,\zeta_2 \in \cCC$ with $\zeta_i|_{[-h,t]} = \xi$
and $\|\zeta_i(s)-\xi(t)\|<\delta$ for all $s\in[t,t+\tau]$ and $i=1,2$, we have
\begin{multline*}
 \underset{s\in [t,t+\tau]}{\esssup}\  \|T(\zeta_1 )(s)-T(\zeta_2) (s)\| \\
 \leq c_0 \sup_{s\in [t,t+\tau]}\|\zeta_1(s)-\zeta_2(s)\|.
\end{multline*}
\item[\textbf{(P3)}]
 $T$ locally maps bounded functions to bounded functions, i.e., for all $\tau>0$ and all $c_1 >0$, there exists $c_2 >0$ such that, for all $\zeta\in C([-h,\tau],\R^n)$,
\[
\sup_{t\in[-h,\tau]}\|\zeta(t)\|\le c_1 ~~\implies~~ \underset{t\in[0,\tau]}{\esssup}\ \|T(\zeta)(t)\| \le c_2.
\]
\end{enumerate}
\end{Def}

We note that an element $T$ of the operator class $\cTT$ is usually the solution operator of a (partial) differential equation describing the internal dynamics of the system. Beyond that, the formulation embraces a large number of processes and effects, such as nonlinear delay elements, backlash and relay hysteresis, and solution operators of infinite-dimensional systems, cf.~\cite{BergIlch21,IlchRyan02b,BergPuch20a}. A practically relevant example where infinite-dimensional internal dynamics appear (and are modelled by an operator $T$) is a moving water tank system considered in~\cite{BergPuch22}. The causality property~(P1) in the class $\cTT$ is physically-motivated and entirely natural. Properties~(P2) and~(P3) are technical conditions required to guarantee the existence of solutions of~\eqref{eq:Sys} under feedback. We stress that property~(P3) is weaker than the respective property required in~\cite{BergIlch21,BergLe18a}, where it essentially needs to hold for ``$\tau=\infty$'' (and hence corresponds to a minimum phase property, cf.\ also Remark~\ref{Rem:LinSys}), while for our purposes a local version suffices.

Next we introduce a sector bound property of~$f\in  C(\R^p\times \R^q\times\R^m,\R^m)$ and~$T\in {\mathbb T}^{rm,q}_h$ as follows.

\begin{enumerate}
\item[\textbf{(P4)}] For all $y^0 \in C^{r-1}([-h,0],\R^m)$ there exist $M_1, \ldots, M_{r+1}\in C(\R_{\ge 0}\times\R^p\times\R^m,\R_{\ge 0})$ such that for all $t\ge 0$, all $(d,v)\in\R^p\times\R^m$ and all $\zeta_1,\ldots,\zeta_r\in C([-h,t],\R^m)$ with $\zeta_i|_{[-h,0]} = (y^0)^{(i-1)}$ for $i=1,\ldots,r$ we have:
    \begin{multline*}
   \hspace{-0.8cm} \|f(d,T(\zeta_1,\ldots,\zeta_r)(t),v)\| \le M_1(t,d,v) \\
    \hspace{-0.4cm}+ M_2(t,d,v) \|\zeta_1|_{[-h,t]}\|_\infty + \ldots + M_{r+1}(t,d,v) \|\zeta_r|_{[-h,t]}\|_\infty
    \end{multline*}
\end{enumerate}

Note that the functions $M_i$ in (P4) depend on the initial history $y^0$ in~\eqref{eq:Sys}.

We are now in the position to define the class of systems to be considered in this paper. We stress that the high-gain property of system~\eqref{eq:Sys} (see property~(P6) below) required in earlier approaches, see e.g.~\cite{BergIlch21}, is not needed here.

\begin{Def}\label{Def:SysClass}
For $m,r\in\N$ we say that system~\eqref{eq:Sys} belongs to the system class $\cN^{m,r}$, written $(d,f,T)\in\cN^{m,r}$, if $d\in L^\infty(\R_{\ge 0},\R^p)$, $f\in C(\R^p\times \R^q\times\R^m,\R^m)$, $T\in {\mathbb T}^{rm,q}_h$ for some $p,q\in\N$, $h\geq 0$ and $(f,T)$ satisfy property~(P4).
\end{Def}

Definition~\ref{Def:SysClass} reflects the spirit of funnel control: The control should work for any system within a class of systems and the class is only described by a set of structural assumptions. If additional information about a specific member of this class is available~-- which we do not assume in this article~-- then the proof of the main result may be tailored to this specific system, if possible.

\begin{Rem}\label{Rem:BIF-nl}
An important subclass of $\cN^{m,r}$ are state-space systems of the form
\begin{equation}\label{eq:sss}
\begin{aligned}
  \dot x(t) &= f(x(t)) + g(x(t)) u(t),\\
  y(t) &= h(x(t)),
\end{aligned}
\end{equation}
where $f:\R^n\to \R^n$, $g:\R^n\to \R^{n\times m}$ and $h:\R^n\to \R^m$ are sufficiently smooth and have the following properties. First of all, the system has relative degree $r\in\N$, which means that $(L_g L_f^k h)(z) = 0$ for all $z\in\R^n$ and all $k=0,\ldots,r-2$ and $(L_g L_f^{r-1} h)(z)$ is invertible for all $z\in\R^n$. Here $(L_f h)(z) := h'(z) f(z)$ denotes the \textit{Lie derivative} of~$h$ along~$f$ at $z\in\R^n$ and we may gradually define $L_f^k h := L_f (L_f^{k-1} h)$ with $L_f^0 h := h$. Furthermore, we set $(L_g h)(z) := [(L_{g_1} h)(z), \ldots, (L_{g_m} h)(z)]$ for $z\in\R^n$, where $g_1(z),\ldots,g_m(z)$ denote the columns of~$g(z)$.

If~\eqref{eq:sss} has relative degree~$r$, then, under the additional technical assumptions of~\cite[Cor.~5.6]{ByrnIsid91a}, there exists a diffeomorphism $\Phi:\R^n\to\R^n$ such that the coordinate transformation $(y(t), \dot y(t),\ldots, y^{(r-1)}(t),\eta(t)) = \Phi(x(t))$ puts the system~\eqref{eq:sss} into Byrnes-Isidori form
\begin{equation}\label{eq:BIF-nl}
\begin{aligned}
  y^{(r)}(t) &= p(y(t), \dot y(t),\ldots, y^{(r-1)}(t),\eta(t))\\
  &\quad  + \gamma(y(t), \dot y(t),\ldots, y^{(r-1)}(t),\eta(t)) u(t),\\
  \dot \eta(t) &= q(y(t), \dot y(t),\ldots, y^{(r-1)}(t),\eta(t)),
\end{aligned}
\end{equation}
where $p:\R^n\to\R^m$, $q:\R^n\to\R^{n-rm}$ and $\gamma = L_g L_f^{r-1} h :\R^n\to \R^{m\times m}$ are continuously differentiable. In order to satisfy~(P4) we assume that~$p$ and~$\gamma$ are globally Lipschitz continuous. The second of equations~\eqref{eq:BIF-nl} is called the \textit{internal dynamics} of~\eqref{eq:sss} and we assume that it satisfies the condition
\begin{multline*}
\hspace{-4mm} \exists\, M\!\in\! C(\R_{\ge 0},\R)\ \forall\,t\!\ge\! 0\ \forall\, \eta^0\!\in\!\R^{n-rm}\ \forall\, \zeta\!\in\! L^\infty_{\loc}([0,t],\R^{rm})\!:\\
 \|\eta(t; \eta^0,\zeta)\| \le M(t)\big( 1 + \|\eta^0\| + \|\zeta|_{[0,t]}\|_\infty\big),
\end{multline*}
where $\eta(\cdot; \eta^0,\zeta)$ denotes the unique solution of $\dot \eta(t) = q(\zeta(t),\eta(t))$ with initial condition $\eta(0) = \eta^0$. Note that under this condition the maximal solution $\eta(\cdot; \eta^0,\zeta)$ can indeed be extended to a global solution. It is now straightforward to check that with the operator
\begin{align*}
  T_{\eta^0}: \zeta \mapsto \begin{pmatrix} \zeta\\ \eta(\cdot; \eta^0,\zeta)  \end{pmatrix}
\end{align*}
and the function $f(d,z,u) = p(z) + \gamma(z) u$ we have that $(0,f,T_{\eta^0})\in \cN^{m,r}$. Therefore, system~\eqref{eq:BIF-nl} belongs to this system class.
\end{Rem}

In contrast to earlier approaches as in~\cite{BergIlch21,BergLe18a}, in this work we consider an additional function $\sat$ in~\eqref{eq:Sys}, which represents an input saturation. If $\sat = {\rm id}_{\R^m}$, then the results from~\cite{BergIlch21,BergLe18a} could be applied. For this reason, we consider a proper input saturation, which has the following, quite general, property.

\begin{enumerate}
\item[\textbf{(P5)}] $\sat\!\in\! C(\R^m,\R^m)$ is bounded and there exists $\theta>0$ such that for all $v\in\R^m$ with $\|v\|\le \theta$ we have $\sat(v) = v$.
\end{enumerate}

We stress that the input saturation function $\sat$ must be known to the controller and it can be viewed as a design parameter, chosen according to the specific requirements of the application at hand. The above property~(P5) allows for a large variety of possible saturations, apart from the standard saturation $\sat_i(v) = v_i$ for $|v_i|\le M$ and $\sat_i(v) = \sgn(v_i) M$ for $|v_i|>M$ for all $i=1,\ldots,m$.

\begin{Rem}\label{Rem:sector}
We like to expound, why the sector bound property~(P4) cannot be dispensed in the presence of (arbitrary) input constraints in general. To this end, consider the following prototype system, where such a linear bound is not satisfied:
\begin{equation}\label{eq:Ex-P4}
    \dot y(t) = y(t)^2 + \sat(v(t)),\quad y(0)=y^0\in\R.
\end{equation}
If $\sat = {\rm id}_\R$, then the classical funnel controller $v(t) = (1-w(t)^2)^{-1}$ with $w(t) = \tfrac{e(t)}{\psi(t)}$ and $e(t)= y(t)-y_{\rm ref}(t)$, for some reference signal~$y_{\rm ref}$ and funnel boundary~$\psi$ with appropriate properties, achieves that the closed-loop differential equation has a global solution~-- see~\cite{BergIlch21} for more details. If
\begin{equation}\label{eq:sat}
 \sat(v) = \begin{cases} v, & |v|\le M,\\ \sgn(v) M,& |v|>M,\end{cases}
\end{equation}
where $M>0$ is some constant, then the application of the above controller (or any other controller) leads to a closed-loop differential equation, which always has a solution with finite escape time, when the saturation level~$M$ is too small.\\
To see this, consider $y^0=1$ and let us assume that the saturation is active with negative sign (a positive control value would only lead to an earlier blow-up), i.e., $\sat(v(t)) = -M$. Then $\dot y(t) \ge y(t)^2 - M$, from which it follows that $y(t) \ge z(t)$ for all $t\in[0,\omega)$ with
\[
    z(t) = \sqrt{M} \frac{\sqrt{M} + 1 + (1-\sqrt{M}) e^{2\sqrt{M} t}}{\sqrt{M}+1 - (1-\sqrt{M}) e^{2\sqrt{M} t}}.
\]
It is straightforward to see that for $M\ge 1$ we have that $z(t)$ is defined for all $t\ge 0$, thus $\omega=\infty$ and hence a global solution exists in this case. However, if $M<1$, then the denominator of $z(t)$ has a zero at
\[
    \omega = \frac{1}{2\sqrt{M}} \ln \left(\tfrac{1+\sqrt{M}}{1-\sqrt{M}}\right),
\]
and hence the solution exhibits a blow-up on the finite interval $[0,\omega)$ in this case. It is clear that, since the maximal possible saturation is already active, no control law would be able to prevent this blow-up; therefore, property~(P4) is mandatory.\\
In fact, blow-up is a typical phenomenon also in higher-dimensional ordinary differential equations of the form $\dot x(t) = p(x(t))$, when~$p$ is a polynomial which (component-wise) involves terms of degree larger than one. The rate of blow-up can even be given in terms of the degree $L = \deg p$: $\|x(t)\| \sim c (T-t)^{-\tfrac{1}{L-1}}$ for $t\to T$; see~\cite[Thm.~4.1]{EliaGing06}.
\end{Rem}

For purposes of comparison we also introduce the system class $\cN^{m,r}_{\rm BIR}$ from~\cite{BergIlch21}, where ``BIR'' stands for the initials of the surnames of the authors of that article. To this end, we recall the high-gain property from~\cite[Def.~1.2]{BergIlch21}.

\begin{enumerate}
\item[\textbf{(P6)}] A function $f\in C(\R^p\times \R^q\times\R^m,\R^m)$ satisfies the \textit{high-gain property}, if there exists $\nu^*\in(0,1)$ such that for all compact $K_p\subset\R^p$ and $K_q\subset \R^q$ the function
\[
    \hspace{-8mm} \chi:\R\!\to\!\R,\ s\!\mapsto\! \min \setdef{\! w^\top f(\delta, z, -sw)}{ \!\!\!\begin{array}{l} (\delta,z)\in K_p\times K_q,\\ w\in\R^m,\\ \nu^* \le \|w\|\le 1\end{array}\!\!\!}
\]
satisfies $\sup_{s\in\R}\chi(s) = \infty$.
\end{enumerate}

With this we may introduce $\cN^{m,r}_{\rm BIR}$ as  follows.

\begin{Def}\label{Def:SysClass-BIR}
For $m,r\in\N$ we say that system~\eqref{eq:Sys} belongs to the system class $\cN^{m,r}_{\rm BIR}$, written $(d,f,T)\in\cN^{m,r}_{\rm BIR}$, if $d\in L^\infty(\R_{\ge 0},\R^p)$, $f\in C(\R^p\times \R^q\times\R^m,\R^m)$, $T\in {\mathbb T}^{rm,q}_h$ for some $p,q\in\N$, $h\geq 0$ and~$T$ satisfies~(P3) for $\tau=\infty$ and~$f$ satisfies property~(P6).
\end{Def}

We stress that in $\cN^{m,r}_{\rm BIR}$ property~(P3) needs to hold with ``$\tau=\infty$'' and hence becomes a bounded-input, bounded-output stability property. For linear systems this property typically corresponds to the minimum phase property as outlined in the following remark.

\begin{Rem}\label{Rem:LinSys} As a consequence of Remark~\ref{Rem:BIF-nl}, the system class $\cN^{m,r}$ in particular contains all linear systems with strict relative degree~$r$, which can be put into Byrnes-Isidori form (cf.~\cite{IlchRyan07,Isid95})
\begin{equation}\label{eq:BIF}
\begin{aligned}
  y^{(r)}(t) &= \sum_{i=1}^r R_i y^{(i-1)}(t) + S\eta(t) + \Gamma u(t),\\
  \dot \eta(t) &= Q\eta(t) + P y(t),
\end{aligned}
\end{equation}
where $R_i, \Gamma\in\R^{m\times m}$, $i=1,\ldots,r$, $S, P^\top \in\R^{m\times q}$ and $Q\in\R^{q\times q}$, with initial conditions $y^{(i-1)}(0) = y_i^0\in\R^m$, $i=1,\ldots,r$, and $\eta(0)=\eta^0\in\R^q$. It is straightforward to check that with the operator
\begin{align*}
  T_{\eta^0}\!:\! (\zeta_1,\ldots,\zeta_r)\!\mapsto\!\! \sum_{i=1}^r R_i \zeta_i \!+\! S e^{Q \cdot} \eta^0 \!+\! S \!\! \int_0^{\cdot} e^{Q(\cdot - s)} P \zeta_1(s) {\rm d}s
\end{align*}
and the function $f(d,z,u) = z + \Gamma u$ we have that $(0,f,T_{\eta^0})\in \cN^{m,r}$.
Note that if all eigenvalues of~$Q$ have negative real part (which means that system~\eqref{eq:BIF} is minimum phase) and~$\Gamma$ is sign-definite (i.e., $\Gamma+\Gamma^\top$ is either positive or negative definite), then we even have $(0,f,T_{\eta^0})\in \cN^{m,r}_{\rm BIR}$, cf.~\cite[Sec.~2.1]{BergIlch21}. Let us emphasize that the latter two requirements are not needed for the class $\cN^{m,r}$: any $Q\in\R^{q\times q}$ and any  $\Gamma\in\R^{m\times m}$ (including $\Gamma = 0$) are allowed within this class. Therefore, the class of linear systems amenable to funnel control by the new controller design presented in this paper is much larger than the class of linear systems considered in earlier works.
\end{Rem}

Remark~\ref{Rem:LinSys} illustrates that the high-gain property~(P6) (typically associated to the concept of ``control direction'') required in~\cite{BergIlch21} for funnel control is not demanded in the system class $\cN^{m,r}$. It is not even required that~$f$ depends on~$u$; however, this would typically mean that the tracking error grows unbounded.

\subsection{Control objective}\label{Ssec:ContrObj}

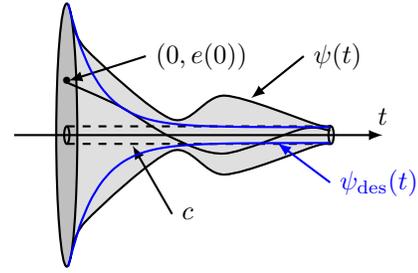
\begin{figure}[h]
  \begin{center}
\begin{tikzpicture}[scale=0.35]
\tikzset{>=latex}
  \filldraw[color=gray!25] plot[smooth] coordinates {(0.15,4.7)(0.7,3.3)(4,0.6)(6,1.5)(9.5,0.4)(10,0.333)(10.01,0.331)(10.041,0.3) (10.041,-0.3)(10.01,-0.331)(10,-0.333)(9.5,-0.4)(6,-1.5)(4,-0.6)(0.7,-3.3)(0.15,-4.7)};
  \draw[thick] plot[smooth] coordinates {(0.15,4.7)(0.7,3.3)(4,0.6)(6,1.5)(9.5,0.4)(10,0.333)(10.01,0.331)(10.041,0.3)};
  \draw[thick] plot[smooth] coordinates {(10.041,-0.3)(10.01,-0.331)(10,-0.333)(9.5,-0.4)(6,-1.5)(4,-0.6)(0.7,-3.3)(0.15,-4.7)};
  \draw[thick,fill=lightgray] (0,0) ellipse (0.4 and 5);
  \draw[thick] (0,0) ellipse (0.1 and 0.333);
  \draw[thick,fill=gray!25] (10.041,0) ellipse (0.1 and 0.333);
  \draw[thick] plot[smooth] coordinates {(0,2)(2,1.1)(4,-0.1)(6,-0.7)(9,0.25)(10,0.15)};
  \draw[thick,->] (-2,0)--(12,0) node[right,above]{\normalsize$t$};
  \draw[thick,dashed](0,0.333)--(10,0.333);
  \draw[thick,dashed](0,-0.333)--(10,-0.333);
  \node [black] at (0,2) {\textbullet};
  \draw[->,thick](4,-3)node[right]{\normalsize$c$}--(2.5,-0.4);
  \draw[->,thick](3,3)node[right]{\normalsize$(0,e(0))$}--(0.07,2.07);
  \draw[->,thick](9,3)node[right]{\normalsize$\psi(t)$}--(7,1.4);
  \draw [color=blue,thick,smooth,domain=0.05:10] plot(\x,{4.7*exp(-0.9*(\x-0.05))+0.3});
  \draw [color=blue,thick,smooth,domain=0.05:10] plot(\x,{-4.7*exp(-0.9*(\x-0.05))-0.3});
  \draw[->,thick,color=blue](10,-1.8)node[right]{\normalsize$\psi_{\rm des}(t)$}--(8,-0.3);
\end{tikzpicture}
\end{center}
 \caption{Error evolution in a funnel $\mathcal F_{\psi}$ with boundary $\psi(t)$ and desired shape $\psi_{\rm des}(t)$.}
 \label{Fig:funnel}
 \end{figure}

The objective is to design a dynamic output derivative feedback strategy such that for any reference signal $y_{\rm ref}\in C^{r}(\R_{\ge 0},\R^m)$ the tracking error $e = y-y_{\rm ref}$ evolves within a performance funnel
\[
    \cF_\psi := \setdef{(t,e)\in\R_{\ge 0}\times\R^m}{\|e\|<\psi(t)},
\]
see Fig.~\ref{Fig:funnel}, which has a desired shape of the form $\psi_{\rm des}(t) = a e^{-bt} +c$ whenever the saturation in~\eqref{eq:IC} is not active, i.e., $\sat(v(t)) = v(t)$, and the actual funnel boundary $\psi(t)$ is allowed to deviate from this shape and become larger when the saturation is active. The specific value of~$\psi(t)$ should be determined by a dynamic part of the control law.

It is usually the hallmark of funnel control that the funnel boundary is prescribed a priori and can be freely chosen by the designer, see e.g.~\cite{BergIlch21,IlchRyan02b,BergLe18a}. Here we do not allow for an arbitrary funnel boundary in order to be able to change its shape by means of a differential equation. However, we allow to prescribe the ``asymptotic shape'' $\psi_{\rm des}(t) = a e^{-bt} +c$ under inactive saturation, that is the positive parameters $a,b,c$ can be chosen as desired.

We also like to note that we do not assume that the reference signal $y_{\rm ref}$ or any of its derivatives is bounded~-- this is also different from classical approaches as mentioned above. In the case of an unbounded reference signal the controller will become saturated at some time and then force the funnel boundary to grow unbounded, thus still guaranteeing the existence of a global solution~-- although with a very bad tracking performance.

\subsection{Organization of the present paper}\label{Ssec:Orga}

The paper is structured as follows. In Section~\ref{Sec:ConStruc}, we introduce a novel funnel controller for systems~\eqref{eq:Sys} under input constraints~\eqref{eq:IC}. Feasibility of the control is proved in the main result in Section~\ref{Sec:Main}: existence of a global solution for systems of class~$\cN^{m,r}$ is shown in Theorem~\ref{Thm:FunCon} and boundedness of this solution (together with convergence of the funnel functions) for sufficiently large saturation level and systems of class~$\cN^{m,r}_{\rm BIR}$ is shown in Theorem~\ref{Thm:FunCon-BIR}. The performance of the funnel controller is compared to that from~\cite{BergIlch21} by an illustrative example in Section~\ref{Sec:Sim}. The paper concludes with Section~\ref{Sec:Concl}.

\section{Funnel control structure}\label{Sec:ConStruc}
%

\begin{figure*}[h!tb]
\centering
\resizebox{\textwidth}{!}{
   \begin{tikzpicture}[very thick,scale=0.7,node distance = 9ex, box/.style={fill=white,rectangle, draw=black}, blackdot/.style={inner sep = 0, minimum size=3pt,shape=circle,fill,draw=black},blackdotsmall/.style={inner sep = 0, minimum size=0.1pt,shape=circle,fill,draw=black},plus/.style={fill=white,circle,inner sep = 0,very thick,draw},metabox/.style={inner sep = 3ex,rectangle,draw,dotted,fill=gray!20!white}]
 \begin{scope}[scale=0.5]
    \node (sys) [box,minimum size=7ex]  {$y^{(r)}(t)= f\big(d(t), T(y,\dot{y},\dots,y^{(r-1)})(t), u(t)\big)$};
    \node [minimum size=0pt, inner sep = 0pt,  below of = sys, yshift=3ex] {System $(d,f,T)\in  \cN^{m,r}$};
    \node(fork1) [minimum size=0pt, inner sep = 0pt,  right of = sys, xshift=45ex] {};
    \node(end1)  [minimum size=0pt, inner sep = 0pt,  right of = fork1, xshift=5ex] {$\big(y,\ldots,y^{(r-1)}\big)$};

   \draw[->] (sys) -- (end1) node[pos=0.4,above] {};

  \node(FC0) [box, below of = fork1,yshift=-8ex,minimum size=7ex] {{$\begin{aligned}
    e_1(t) &= e(t) = y(t) - y_{\rm ref}(t)\\
    e_{i+1}(t) &= e^{(i)}(t) + k_i(t) e_i(t)\end{aligned}$}};
    \node(FC1) [box, below of = sys,yshift=-8ex,minimum size=7ex] {$\begin{aligned}\dot \psi_i(t) &= p_i \psi_{i+1}(t) - \alpha_i \psi_i(t) + \beta_i - p_i \tfrac{\beta_{i+1}}{\alpha_{i+1}},\\
    \dot \psi_r(t) &= -\alpha_r \psi_r(t) + \beta_r + \psi_r(t) \tfrac{\kappa(v(t))}{\|e_r(t)\|}\end{aligned}$};
    \draw[->] (fork1) -- (FC0) {};
    \node(FC2) [box, below of = FC1,yshift=-5ex,minimum size=7ex] {$v(t) =N\big(k_r(t)\big) e_r(t)$};
     \node(sat) [box, left of = FC1,xshift=-25ex,minimum size=7ex] {$u(t) = \sat(v(t))$};
   \draw[->] (FC2) -- (FC1) node[midway,right] {$v$};
   \node(e0) [minimum size=0pt, inner sep = 0pt,  below of = FC0, xshift=-7ex, yshift=5ex] {};
   \node(e1) [minimum size=0pt, inner sep = 0pt,  below of = FC0, xshift=7ex, yshift=5ex] {};
   \draw[->] (e0) |- (FC2) node[pos=0.7,above] {$e_r$};
   \node(phi) [minimum size=0pt, inner sep = 0pt,  below of = e1, yshift=0ex] {$\big(y_{\rm ref},\ldots,y_{\rm ref}^{(r-1)}\big)$};
   \draw[->] (phi) -- (e1) node[midway,right] {};
   \draw[->] (FC2) -| (sat)  node[pos=0.7,right] {$v$};
   \draw[->] (sat) |- (sys) node[pos=0.7,above] {$u$};
   \node(f0) [minimum size=0pt, inner sep = 0pt,  left of = FC0, xshift=-5.9ex, yshift=2ex] {};
   \node(f1) [minimum size=0pt, inner sep = 0pt,  left of = FC0, xshift=-5.9ex, yshift=-2ex] {};
   \node(f2) [minimum size=0pt, inner sep = 0pt,  left of = f0, xshift=-9.1ex, yshift=0ex] {};
   \node(f3) [minimum size=0pt, inner sep = 0pt,  left of = f1, xshift=-9.1ex, yshift=0ex] {};
   \draw[<-] (f3) -- (f1) node[midway,below] {$e_r$};
   \draw[<-] (f0) -- (f2) node[midway,above] {$\psi_1,\ldots,\psi_r$};

   \node [minimum size=0pt, inner sep = 0pt,  below of = FC2, yshift=1ex, xshift=-12ex] {Funnel controller~\eqref{eq:ICFC}};
\end{scope}
\begin{pgfonlayer}{background}
      \fill[lightgray!20] (-5.3,-2.2) rectangle (16.1,-9.2);
      \draw[dotted] (-5.3,-2.2) -- (16.1,-2.2) -- (16.1,-9.2) -- (-5.3,-9.2) -- (-5.3,-2.2);
  \end{pgfonlayer}
  \end{tikzpicture}
}
\caption{Construction of the funnel controller~\eqref{eq:ICFC} and its internal feedback loops.}
\label{Fig:funnel-controller}
\end{figure*}
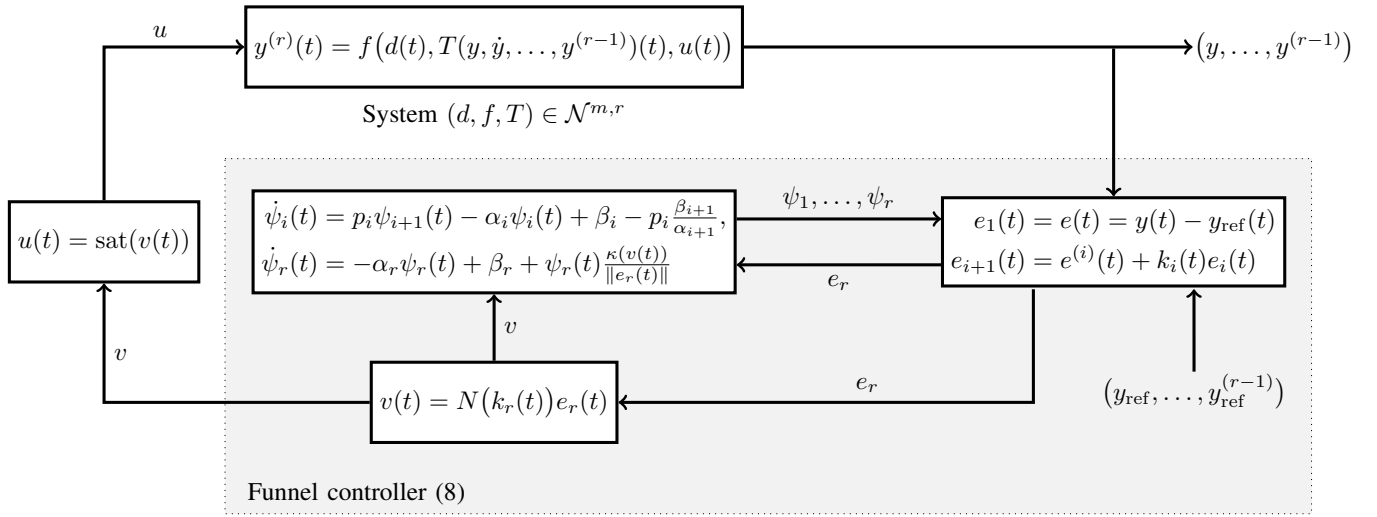

We introduce the following input-constrained funnel controller for systems~\eqref{eq:Sys},~\eqref{eq:IC}.

\begin{equation}\label{eq:ICFC}
\boxed{
\begin{aligned}
    e_1(t) &= e(t) = y(t) - y_{\rm ref}(t),\\
    e_{i+1}(t) &= e^{(i)}(t) + k_i(t) e_i(t),\quad i=1,\ldots,r-1,\\
    k_i(t) &= \left(1- \frac{\|e_i(t)\|^2}{\psi_i(t)^2}\right)^{-1},\quad i=1,\ldots,r,\\
    \dot \psi_i(t) &= p_i \psi_{i+1}(t) - \alpha_i \psi_i(t) + \beta_i - p_i \frac{\beta_{i+1}}{\alpha_{i+1}},\\
    \psi_i(0)&= \psi_i^0,\qquad\qquad\qquad i=1,\ldots,r-1,\\
    \dot \psi_r(t) &= -\alpha_r \psi_r(t) + \beta_r + \psi_r(t) \frac{\kappa(v(t))}{\|e_r(t)\|},\\
    \psi_r(0) &= \psi_r^0,\\
    \kappa(v(t)) &= \|v(t)-\sat(v(t))\|,\\
    v(t) &= N\big(k_r(t)\big) e_r(t)
\end{aligned}
}
\end{equation}
with the controller design parameters
\begin{equation}\label{eq:FC-param}
\boxed{
\begin{aligned}
    &\alpha_1\!>\!\alpha_2\!>\!\ldots\!>\!\alpha_r\!>\!0,\ p_i\!>\!1\ \text{ for } i=1,\ldots,r\!-\!1, \\
    &\beta_i\!>\!0,\  \psi_i^0\!>\!\frac{\beta_i}{\alpha_i}\ \text{ for } i=1,\ldots,r,\\
    &N\in C(\R_{\ge 0},\R)\ \text{ a surjection}.
\end{aligned}
}
\end{equation}
Furthermore, with reference to Fig.~\ref{Fig:funnel-controller}, in~\eqref{eq:ICFC} we assume that the instantaneous values of the tracking error $e(t)$ and its derivatives $\dot e(t),\ldots, e^{(r-1)}(t)$ are available for feedback, thus~\eqref{eq:ICFC} is a dynamic error derivative feedback controller.

\captionsetup[subfloat]{labelformat=empty}

The first three equations of the controller~\eqref{eq:ICFC} are basically a combination of the two designs from~\cite{BergIlch21,BergLe18a}, appended by the dynamics for the funnel boundaries in the subsequent three equations. This contrasts classical funnel control approaches, where the performance funnels are always prescribed a priori. Here, they are determined by a dynamical system, which is influenced by the input and an auxiliary error variable. Since the funnel functions are then used to determine these quantities in turn, a feedback structure arises (depicted in Fig.~\ref{Fig:funnel-controller}), for which we seek to prove existence of global solutions.

The surjective function~$N$ in~\eqref{eq:FC-param} serves the purpose of accommodating for possibly unknown control directions. With its help the controller is able to ``probe'' for the appropriate sign of the control signal. A typical choice for~$N$ would be $N(s) = s \sin s$. For more details see also~\cite[Rem.~1.8]{BergIlch21}.

The distinguishing feature of the novel control design~\eqref{eq:ICFC} is that it is feasible under arbitrary input constraints~\eqref{eq:IC}. The controller~\eqref{eq:ICFC} always guarantees the evolution of the tracking error within a performance funnel, whose boundary is determined by a dynamic part of the controller as mentioned above. The term $\kappa(v(t))$ in the differential equation for $\psi_r$ determines whether the saturation is active (i.e., $\kappa(v(t))\neq 0$) or inactive (i.e., $\kappa(v(t))= 0$). If the saturation is inactive, then $\psi_r(t) = \left(\psi_r^0-\tfrac{\beta_r}{\alpha_r}\right) e^{-\alpha_r t} + \tfrac{\beta_r}{\alpha_r}$; if the saturation is active, then $\kappa(v(t))$ provides a positive contribution to $\dot \psi_r$ and hence widens the funnel~-- the larger the deviation between $v(t)$ and $\sat(v(t))$, the larger the widening effect. If $\psi_r$ thus deviates from it's desired shape it will contribute a larger positive part to $\dot \psi_{r-1}$ and hence force the funnel boundary $\psi_{r-1}$ to widen. This effect propagates through the dynamics of the funnel boundaries back to $\psi_1$. After a period of active saturation, the boundaries recover to their prescribed shape exponentially fast.

We emphasize that the controller~\eqref{eq:ICFC} introduces several possible singularities in the closed-loop differential equation, via the gain functions~$k_i(t)$ (when $\|e_i(t)\| = \psi_i(t)$) and via the last term in the expression for $\dot \psi_r(t)$ (when $\|e_r(t)\|=0$). In order to prove the existence of a global solution, it must be ensured that $\|e_i(t)\| \le \eps_i \psi_i(t)$ for some $\eps_i\in(0,1)$ and that $\kappa(v(t)) = 0$ whenever $\|e_r(t)\| < \delta$ for some $\delta>0$. Furthermore, compared to classical funnel control approaches as in~\cite{BergIlch21,BergLe18a}, the funnel boundaries~$\psi_i$ are not prescribed here, and in particular it is not known a priori that they are bounded. Hence, solutions may potentially get unbounded in finite time, i.e., exhibit a blow-up. Therefore, the feasibility proof of the control design is a highly nontrivial task.

\section{Funnel control -- main results}\label{Sec:Main}
%

In this section we show that the application of the funnel controller~\eqref{eq:ICFC} to a system~\eqref{eq:Sys} under input constraints~\eqref{eq:IC} leads to a closed-loop initial-value problem which has a global solution.  By a solution of~\eqref{eq:Sys},~\eqref{eq:IC},~\eqref{eq:ICFC} on $[-h,\omega)$ we mean a tuple of functions $(y,\psi_1,\ldots,\psi_r)\in C^{r-1}([-h,\omega),\R^m)\times C([-h,\omega),\R)^r$ with $\omega\in (0,\infty]$, which satisfies $y|_{[-h,0]} = y^0$, $\psi_i(0)=\psi_i^0$ for all $i=1,\ldots,r$ and $(y^{(r-1)},\psi_1,\ldots,\psi_r)|_{[0,\omega)}$ is locally absolutely continuous and
satisfies the differential equations in~\eqref{eq:Sys} and~\eqref{eq:ICFC} with $u$ defined by~\eqref{eq:IC},~\eqref{eq:ICFC} for almost all $t\in[0,\omega)$;
$(y,\psi_1,\ldots,\psi_r)$ is called maximal, if it has no right extension that is also a solution.

Next we present the main result of the present paper.

\begin{Thm}\label{Thm:FunCon}
Consider a system~\eqref{eq:Sys} with $(d,f,T)\in\cN^{m,r}$ for $m,r\in\N$, under input saturation~\eqref{eq:IC} with saturation function $\sat$ that satisfies~(P5). Let $y^0\in C^{r-1}([-h,0],\R^m)$ be the initial trajectory, $y_{\rm ref}\in C^r(\R_{\ge 0},\R^m)$ the reference signal and choose funnel control design parameters as in~\eqref{eq:FC-param}. Set $e=y-y_{\rm ref}$ and assume that the instantaneous values $e(t), \dot e(t),\ldots, e^{(r-1)}(t)$ are available for feedback and satisfy, using the variables $e_1,\ldots,e_r$ defined in~\eqref{eq:ICFC}, that
\begin{equation}\label{eq:InitCond}
    \forall\, i=1,\ldots,r:\ \|e_i(0)\| < \psi_i^0.
\end{equation}
Then the funnel controller~\eqref{eq:ICFC} applied to~\eqref{eq:Sys},~\eqref{eq:IC} yields an initial-value problem which has
a solution, every solution can be maximally extended and every maximal solution $(y,\psi_1,\ldots,\psi_r):[-h,\omega)\to\R^{m+r}$, $\omega\in(0,\infty]$, has the following properties:
\begin{enumerate}
  \item global existence: $\omega = \infty$;
  \item the functions $e_1,\ldots,e_r$ evolve in their respective performance funnels in the sense:
  \begin{multline*}
    \forall\, i=1,\ldots,r-1\ \exists\, \eps_i\in(0,1)\ \forall\, t\ge 0:\\
     \|e_i(t)\| \le  \eps_i \psi_i(t) \quad \text{and}\quad \|e_r(t)\| < \psi_r(t);
  \end{multline*}
  \item if the saturation is not active on some interval $[t_0,t_1)\subseteq\R_{\ge 0}$ with $t_1\in (t_0,\infty]$, i.e., $v(t) = \sat(v(t))$ for all $t\in [t_0,t_1)$, then the performance funnels exponentially recover to their prescribed shape:
      \begin{multline*}
        \forall\, i=1,\ldots,r\ \forall\, t\in [t_0,t_1):\\ \psi_i(t) \le \frac{\beta_i}{\alpha_i} + \sum_{j=i}^r \mu_j(t_0) \nu_{ij} e^{-\alpha_j (t-t_0)},
      \end{multline*}
      where $\mu_i(t_0) := \psi_i(t_0) -  \tfrac{\beta_i}{\alpha_i}$, $\nu_{ii}:=1$ and $\nu_{ij} := \prod_{k=i}^{j-1} \tfrac{p_k}{\alpha_k-\alpha_j}$ for $i=1,\ldots,r$ and $j=i+1,\ldots,r$.
\end{enumerate}
\end{Thm}

The proof is relegated to Appendix~\ref{app:proof_main}.

We stress that although Theorem~\ref{Thm:FunCon} provides the existence of a global solution of the closed-loop system, it cannot be concluded that the funnel boundaries $\psi_1,\ldots,\psi_r$ are bounded in general. However, statement~(iii) provides that \textit{a posteriori} the funnel boundaries recover to their prescribed shape on any interval where the saturation is not active; in particular, if $t_1=\infty$, then they are bounded.

Nevertheless, it is possible to show global boundedness of $\psi_1,\ldots,\psi_r$ for sufficiently large saturation level, i.e., $\sat(v)=v$ for all $v\in\R^m$ with $\|v\|\le M$ and $M>0$ sufficiently large. For this we require additional assumptions, i.e., a bounded reference signal with bounded derivatives and the system class $\cN^{m,r}_{\rm BIR}$ from Definition~\ref{Def:SysClass-BIR}.

\begin{Thm}\label{Thm:FunCon-BIR}
Consider a system~\eqref{eq:Sys} with $(d,f,T)\in\cN^{m,r}_{\rm BIR}$ for $m,r\in\N$. Choose funnel control design parameters as in~\eqref{eq:FC-param}, $\eps\in (0,1)$ and $K>0$. Then there exists $M>0$ (depending on $\eps$ and $K$) such that
\begin{itemize}
  \item for all saturation functions $\sat$ which satisfy~(P5) with $\theta=M$,
  \item for all $y^0\in C^{r-1}([-h,0],\R^m)$ with $\|e_i(0)\|\le \eps \psi_i^0$, $i=1,\ldots,r$, and
  \item for all $y_{\rm ref}\in W^{r,\infty}(\R_{\ge 0},\R^m)$ with $\|y_{\rm ref}^{(i)}\|_\infty \le K$, $i=0,\ldots,r$,
\end{itemize}
there exists a solution~$(y,\psi_1,\ldots,\psi_r):[-h,\omega)\to\R^{m+r}$, $\omega\in(0,\infty]$, of~\eqref{eq:Sys},~\eqref{eq:IC},~\eqref{eq:ICFC} which can be maximally extended to a global solution (i.e., $\omega=\infty$) that satisfies
\begin{enumerate}
  \item $y\in W^{r,\infty}([-h,\infty),\R^m)$ and, for all $i=1,\ldots,r$, $\psi_i\in L^\infty(\R_{\ge 0},\R)$ with
  \[
    \limsup_{t\to\infty} \psi_i(t) \le \frac{\beta_i}{\alpha_i};
  \]
  \item $k_i\in L^\infty(\R_{\ge 0},\R)$ for $i=1,\ldots,r$ and $v\in L^\infty(\R_{\ge 0},\R^m)$ with $\|v(t)\|< M$ for all $t\ge 0$, for the quantities defined in~\eqref{eq:ICFC};
  \item in particular, for $i=1,\ldots,r$ there exists $\eps_i\in (0,1)$ such that $\|e_i(t)\|\le \eps_i \psi_i(t)$ for all $t\ge 0$.
\end{enumerate}
\end{Thm}

The proof is relegated to Appendix~\ref{app:proof_Cor}.

We like to point out that by statement~(ii) of Theorem~\ref{Thm:FunCon-BIR} the saturation is never active and hence the funnel boundaries do not adjust themselves, but follow their desired shape. Also note that the proof of Theorem~\ref{Thm:FunCon-BIR} is constructive, the saturation level~$M$ (depending on~$\eps$ and~$K$) is provided explicitly. We further stress that Theorem~\ref{Thm:FunCon-BIR} provides an explicit relation between the initial values and the saturation level~$M$. The initial values $y(0),\dot y(0),\ldots, y^{(r-1)}(0)$ are essentially confined to a bounded set, the size of which is quantified by $\eps\in(0,1)$, for which the relations $\|e_i(0)\|\le \eps \psi_i^0$ hold for $i=1,\ldots,r$. If $\eps$ is made smaller, allowing only a smaller set of initial values, then it is also possible to choose a smaller saturation level~$M$ in general (although $M\ge M^*>0$ even for $\eps\to 0$ and $K\to 0$, where $M^*$ depends on the system and controller parameters).

With the same proof technique as in Appendix~\ref{app:proof_Cor} it is possible to show that there exists an invariant set for the coordinates $e_i(t)/\psi_i(t)$ within which the saturation is never active: There are $\eps_1,\ldots,\eps_r$ such that whenever $\|e_i(0)\| \le \eps_i \psi_i(0)$ we have $\|e_i(t)\| \le \eps_i \psi_i(t)$ and $\|v(t)\| \le M$ for all $t\ge 0$.

\begin{Cor}\label{Cor:FunCon-Inv}
Consider a system~\eqref{eq:Sys} with $(d,f,T)\in\cN^{m,r}_{\rm BIR}$ for $m,r\in\N$. Choose funnel control design parameters as in~\eqref{eq:FC-param}. Then for all $K>0$ there exist $\hat \eps_1,\ldots,\hat \eps_r\in (0,1)$ such that for all $\eps_i\in [\hat \eps_i, 1)$, $i=1,\ldots,r$ there exists $M>0$ such that
\begin{itemize}
  \item for all saturation functions $\sat$ which satisfy~(P5) with $\theta=M$,
  \item for all $y^0\in C^{r-1}([-h,0],\R^m)$ with $\|e_i(0)\|\le \eps_i \psi_i^0$, $i=1,\ldots,r$, and
  \item for all $y_{\rm ref}\in W^{r,\infty}(\R_{\ge 0},\R^m)$ with $\|y_{\rm ref}^{(i)}\|_\infty \le K$, $i=0,\ldots,r$,
\end{itemize}
there exists a solution~$(y,\psi_1,\ldots,\psi_r):[-h,\omega)\to\R^{m+r}$, $\omega\in(0,\infty]$, of~\eqref{eq:Sys},~\eqref{eq:IC},~\eqref{eq:ICFC} which can be maximally extended to a global solution (i.e., $\omega=\infty$) that satisfies~(i) and~(ii) as in Theorem~\ref{Thm:FunCon-BIR} and additionally we have
\begin{enumerate}
  \item[(iii)] $\|e_i(t)\|\le \eps_i \psi_i(t)$ for all $t\ge 0$ and all $i=1,\ldots,r$.
\end{enumerate}
\end{Cor}

The proof is relegated to Appendix~\ref{app:proof_Cor_2}.

The purpose of Corollary~\ref{Cor:FunCon-Inv} is to obtain recursive feasibility for the extension of Funnel MPC 
to higher relative degree systems and is used for the proof of~\cite[Thm.~2.3]{BergDenn22}.

\section{Simulations}\label{Sec:Sim}

\begin{figure}[htp]
    \begin{center}
    \includegraphics[trim=2cm 4cm 5cm 15cm,clip=true,width=6.5cm]{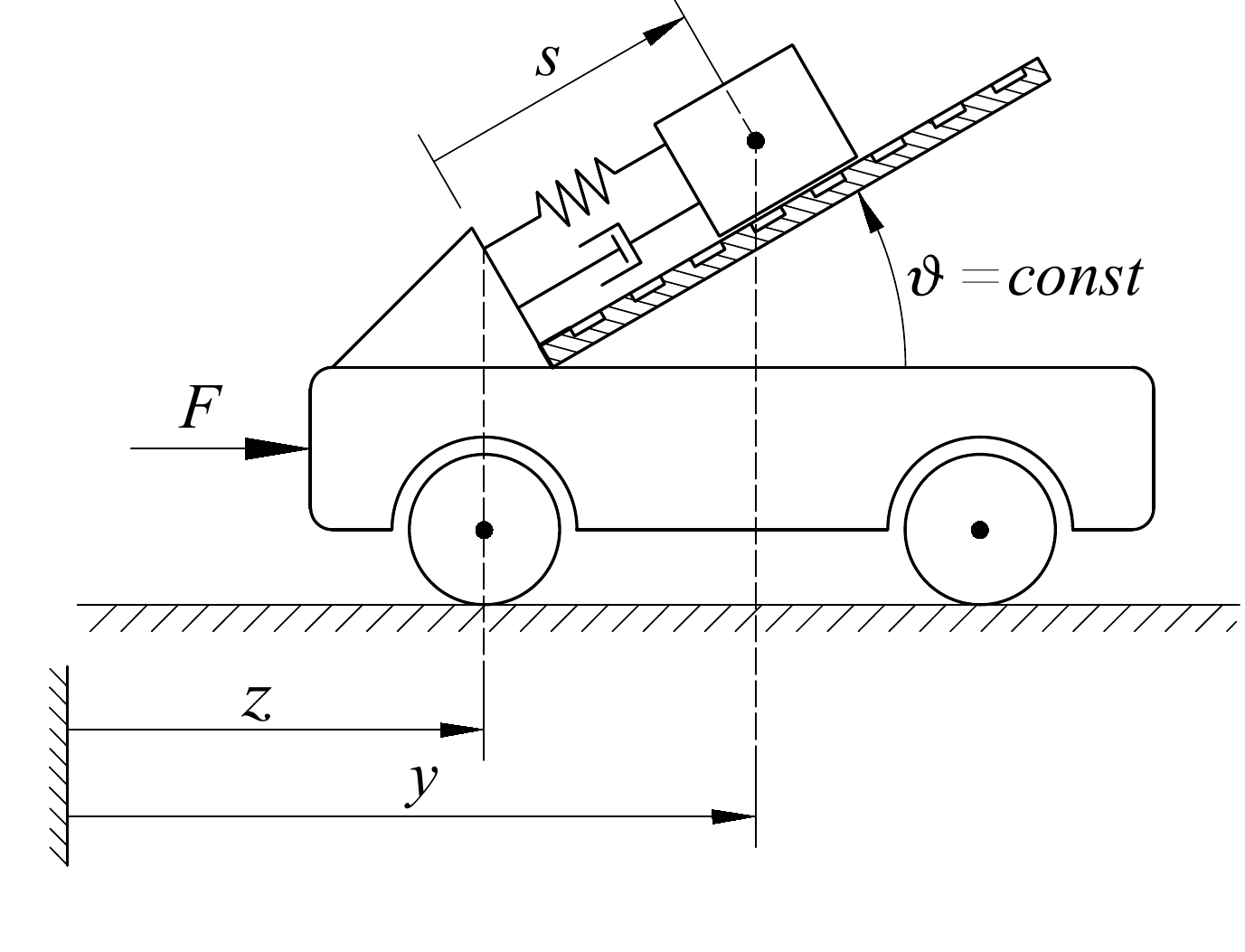}
    \end{center}
    \vspace*{-3mm}
    \caption{Mass-on-car system.}
    \label{Mass.on.car}
    \end{figure}

We compare the controller~\eqref{eq:ICFC} to the controllers presented in~\cite{BergIlch21,BergLe18a} and, to this end, consider the mass-on-car system example presented therein, which is from~\cite{SeifBlaj13}. As illustrated in Fig.~\ref{Mass.on.car}, the mass~$m_2$ (in~\si{\kilo\gram})  moves on a ramp inclined by the angle~{$\vartheta \in [0,\frac{\pi}{2})$} (in \si{\radian}) and is mounted on a car with mass~$m_1$ (in \si{\kilo\gram}). We assume that the control input is the force~$u=F$ (in \si{\newton}) acting on the car. The equations of motion for the system are given by
\begin{equation}\label{mass.on.car.equ}
\begin{bmatrix}
m_1+m_2&m_2\cos \vartheta\\
m_2\cos \vartheta&m_2
\end{bmatrix} \begin{pmatrix} \ddot{{ z}}(t)\\ \ddot{s}(t) \end{pmatrix} +\begin{pmatrix}
0\\
ks(t)\!+\!d\dot{s}(t)
\end{pmatrix}\!=\!\begin{pmatrix}
u(t)\\
0
\end{pmatrix},
\end{equation}
where $t$ is the current time (in \si{\second}),  $z$ (in \si{\metre}) is the horizontal car position and~$s$ (in \si{\metre}) the relative position of the mass on the ramp. The constants~$k  >0$ (in \si{\newton\per{\metre}}), $d >0$ (in \si{\newton\second\per{\metre}}) are the coefficients of the spring and damper, resp. The output
 $y$ (in \si{\metre}) is the horizontal position of the mass on the ramp,
\[
y(t)={ z}(t)+s(t)\cos \vartheta.
\]
It is easy to see that~\eqref{mass.on.car.equ} with output~$y$ can be transformed into the form~\eqref{eq:BIF} and hence, as outlined in Remark~\ref{Rem:LinSys}, belongs to the class $\cN^{1,r}$ with $r=2$ if $\vartheta\in (0,\frac{\pi}{2})$, and $r=3$ if $\vartheta=0$. As derived in~\cite[Sec.~3.1]{BergIlch21}, it even belongs to $\cN^{1,r}\cap \cN^{1,r}_{\rm BIR}$ in both cases.

For the simulation, we choose the parameters
$m_1=4$, $m_2=1$, $k=2$, $d=1$,
the initial values $z(0)= s(0) = 0$, $\dot{z}(0) = \dot s(0)= 0$ and the reference signal $y_{\rm ref} \colon t\mapsto \cos t$. All simulations are MATLAB generated (solver: {\tt ode45}, rel.\ tol.: $10^{-10}$, abs.\ tol.: $10^{-8}$) and over the time interval $[0,15]$. We consider two cases.

\vspace{2mm}
\noindent
{\bf Case 1:} $\vartheta = \frac{\pi}{4}$. Then system~\eqref{mass.on.car.equ} belongs to $\cN^{1,2}$. For the controller~\eqref{eq:ICFC} we choose the controller design parameters from~\eqref{eq:FC-param} as
\begin{align*}
    \alpha_1 &= 1.5,\ \alpha_2 = 0.9\cdot \alpha_1,\ \beta_1 = 0.15,\ \beta_2 = 0.5\cdot \alpha_2, \\
    p_1 &=  1.1,\ \psi_1^0 = 4.1,\ \psi_2^0 = 2
\end{align*}
and $N(s) = -s^2 \cos s$. The saturation function in~\eqref{eq:IC} is chosen as in~\eqref{eq:sat} with $M=10$.

The controller from~\cite{BergIlch21} takes the form
\begin{equation}\label{eq:FC-BIR}
\begin{aligned}
 w(t) &= \varphi(t) \dot e(t) + \alpha\big(\varphi(t)^2 e(t)^2\big)\, \varphi(t) e(t),\\
u(t) &= N\big(\alpha\big(w(t)^2\big)\big)\, w(t),
\end{aligned}
\end{equation}
where~$\alpha(s) = 1/(1-s)$ for~$s\in[0,1)$. We choose~$\varphi(t) = (4 e^{-3t/2} + 0.1)^{-1}$ for $t\geq 0$, so that $\ddt\big(1/\varphi(t)\big) = -\alpha_1 /\varphi(t) + \beta_1$.

\captionsetup[subfloat]{labelformat=empty}
\begin{figure}[h!tb]
  \centering
  \subfloat[Fig.~\ref{fig:sim-rd2}a: Performance funnels and tracking errors]
{
\centering
\hspace{-2mm}
  \includegraphics[width=0.5\textwidth]{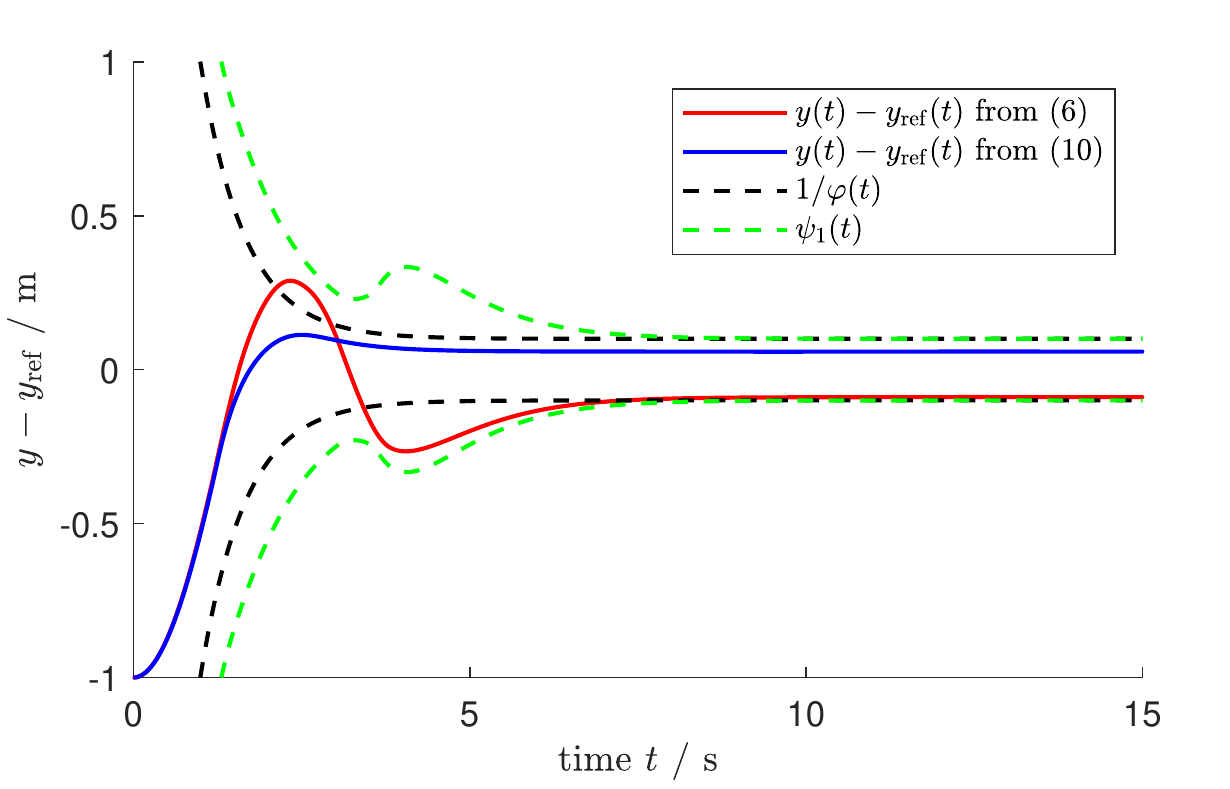}
\label{fig:sim-e}
}\\
\subfloat[Fig.~\ref{fig:sim-rd2}b: Input functions]
{
\centering
\hspace{-5mm}
  \includegraphics[width=0.52\textwidth]{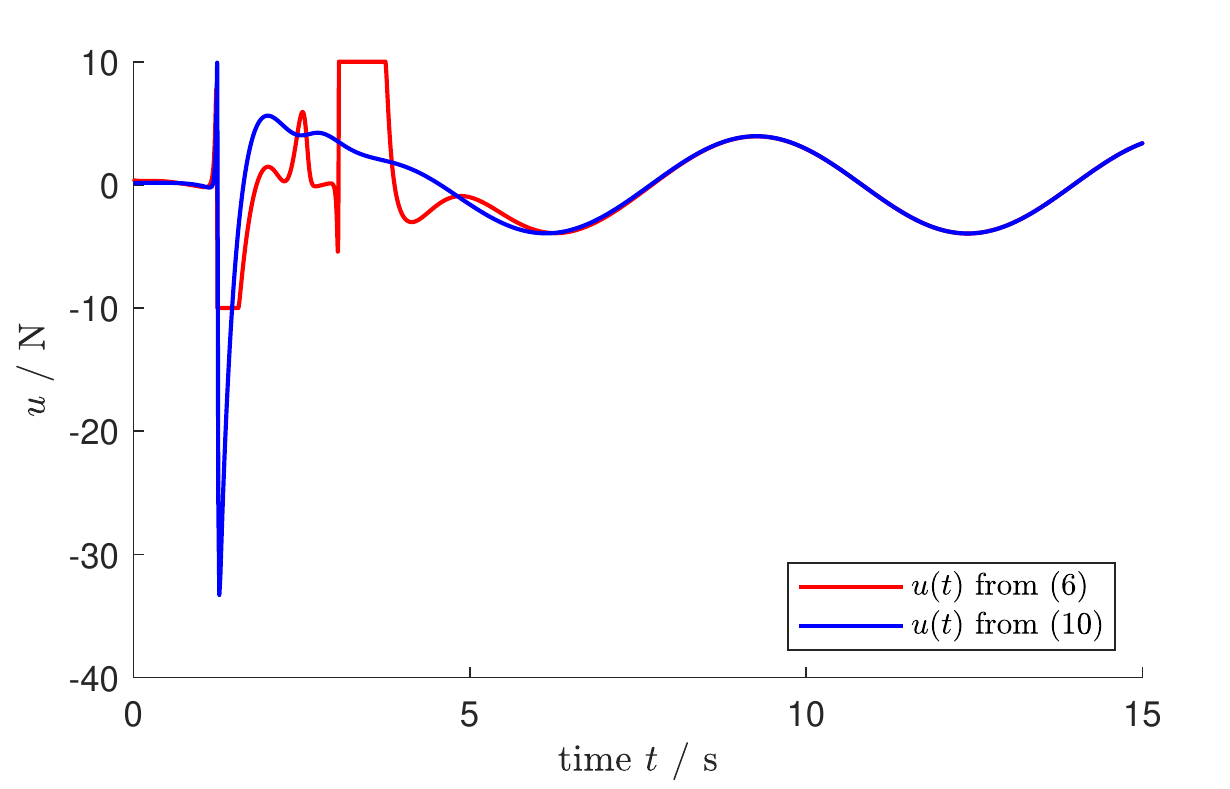}
\label{fig:sim-u}
}
\caption{Simulation, under controllers~\eqref{eq:ICFC} and~\eqref{eq:FC-BIR}, of system~\eqref{mass.on.car.equ} with $\vartheta=\frac{\pi}{4}$.}
\label{fig:sim-rd2}
\end{figure}

The application of the controllers~\eqref{eq:ICFC} and~\eqref{eq:FC-BIR} to~\eqref{mass.on.car.equ} is depicted in Fig.~\ref{fig:sim-rd2}. The corresponding tracking errors and funnel boundaries are shown in
Fig.~\ref{fig:sim-e}, while Fig.~\ref{fig:sim-u} shows the respective input functions. It is evident that the performance of both controllers is comparable (if not identical) whenever the saturation is not active. When the saturation is active the tracking error generated by the new controller~\eqref{eq:ICFC} leaves the performance funnel $\cF_{1/\varphi}$, but stays within the widened funnel $\cF_{\psi_1}$. It can be seen that the enforced widening of the funnel is not very significant, but the saturated control signal of~\eqref{eq:ICFC} avoids the undesirable peak of~\eqref{eq:FC-BIR} (which is a common phenomenon in earlier funnel control approaches). Furthermore, it is important to note that after the saturation was last active on (approximately) $[3,3.8]$ the funnel boundary~$\psi_1$ exponentially converges to $1/\varphi$, thus retaining the tracking error within the desired region again.
\\[2mm]
\noindent
{\bf Case 2:} $\vartheta=0$. Then system~\eqref{mass.on.car.equ} belongs to $\cN^{1,3}$. For the controller~\eqref{eq:ICFC} we choose the parameters from~\eqref{eq:FC-param} as
\begin{align*}
    \alpha_1 &= 1.5,\ \alpha_2 = 0.9\cdot \alpha_1,\ \alpha_3 = 0.9\cdot \alpha_2,\\
     \beta_1& = 0.1,\ \beta_2 = 0.5\cdot \alpha_2,\ \beta_3 = 0.5\cdot \alpha_3\\
    p_1 &= p_2 =  1.1,\ \psi_1^0 = 3.1,\ \psi_2^0 = \psi_3^0 = 1.6
\end{align*}
and again $N(s) = -s^2 \cos s$. The saturation function is chosen as in~\eqref{eq:sat} with $M=8$.

The funnel controller from~\cite{BergIlch21} with $r=3$ takes the form
\begin{equation}\label{eq:FC-BIR-3}
\begin{aligned}
  w(t) &= \varphi(t) \ddot e(t) + {\gamma}\big(\varphi(t) \dot e(t) + {\gamma}\big(\varphi(t) e(t)\big)\big),\\
  u(t) &= N\big(\alpha\big(w(t)^2\big)\big)\, w(t),
\end{aligned}
\end{equation}
where $\gamma(s) = \alpha(s^2) s$ for $s\in(-1,1)$, and we choose $\varphi(t) = (3 e^{-t} + 0.1)^{-1}$ for $t\ge 0$ and
compare the controller~\eqref{eq:ICFC} with~\eqref{eq:FC-BIR-3}.

\captionsetup[subfloat]{labelformat=empty}
\begin{figure}[h!tb]
  \centering
  \subfloat[Fig.~\ref{fig:sim-rd3}a: Performance funnels and tracking errors]
{
\centering
\hspace{-2mm}
  \includegraphics[width=0.5\textwidth]{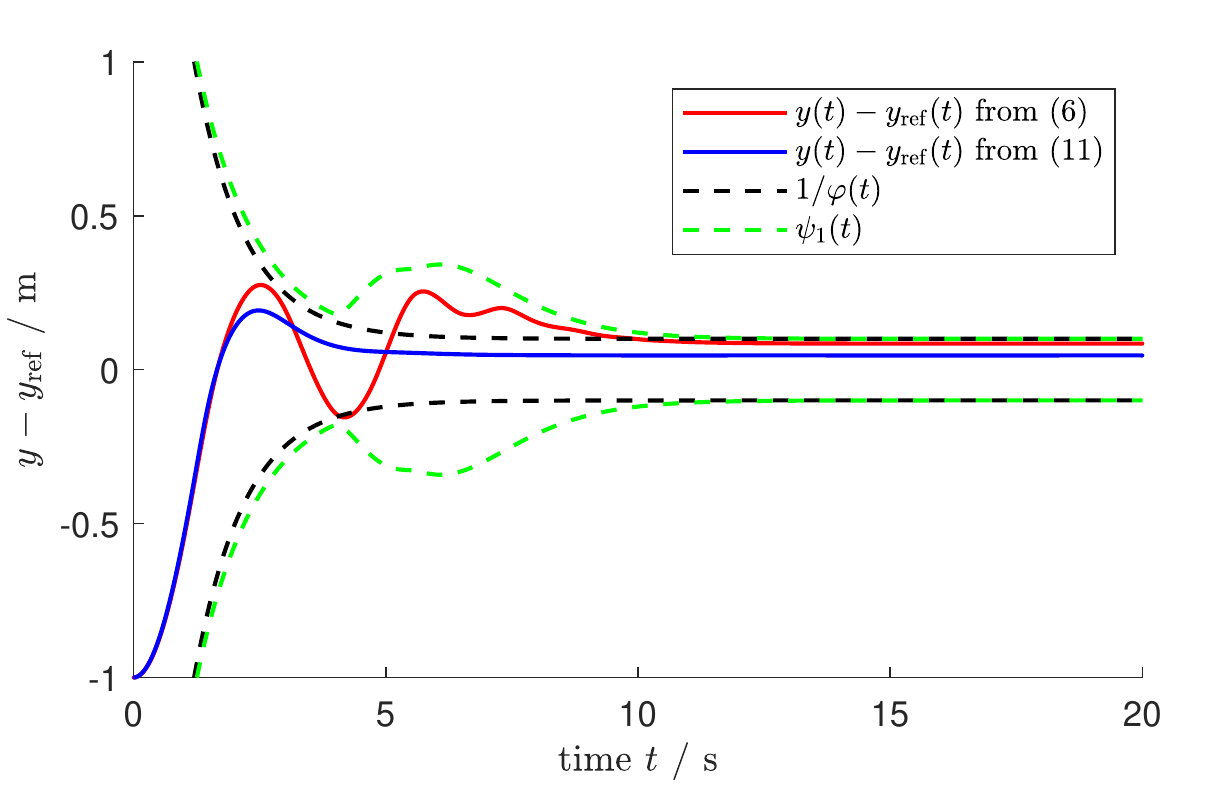}
\label{fig:sim2-e}
}\\
\subfloat[Fig.~\ref{fig:sim-rd3}b: Input functions]
{
\centering
\hspace{-5mm}
  \includegraphics[width=0.52\textwidth]{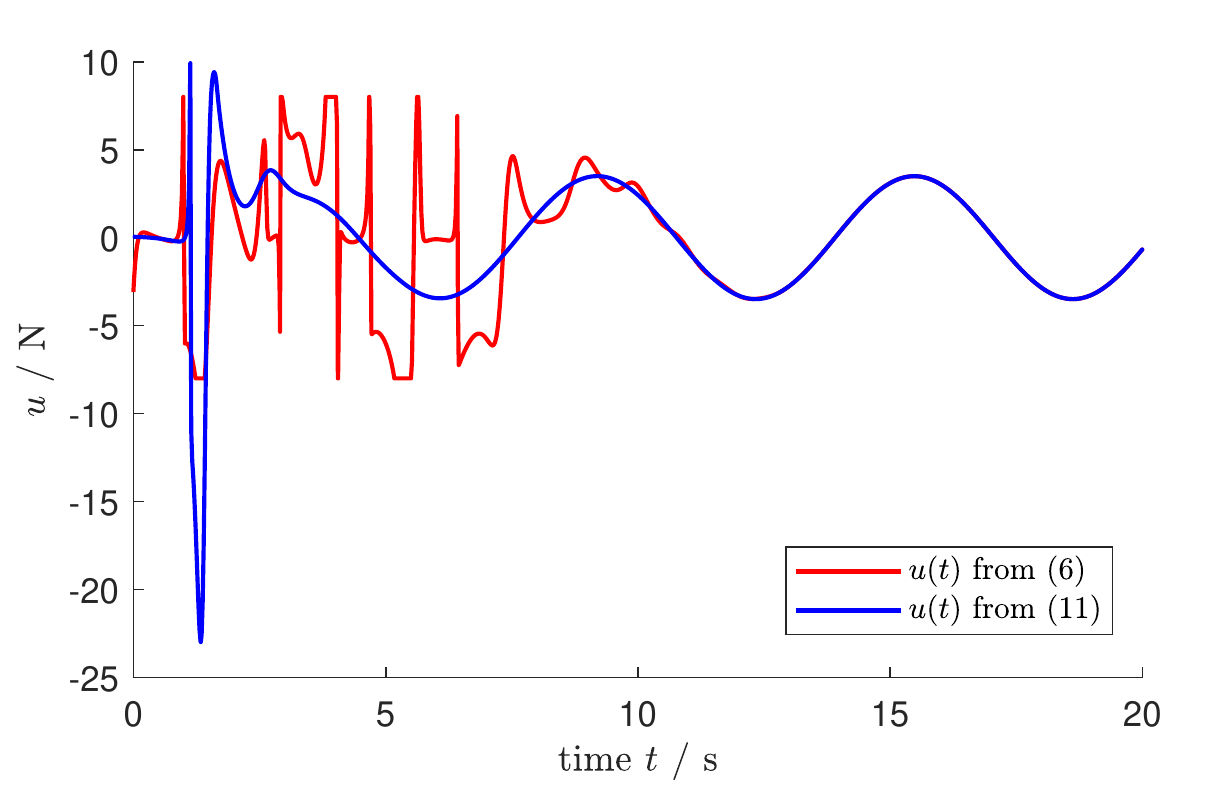}
\label{fig:sim2-u}
}
\caption{Simulation, under controllers~\eqref{eq:ICFC} and~\eqref{eq:FC-BIR-3}, of system~\eqref{mass.on.car.equ} with $\vartheta=0$.}
\label{fig:sim-rd3}
\end{figure}

The simulation depicted in Fig.~\ref{fig:sim-rd3} shows that the controller~\eqref{eq:FC-BIR-3} generates unnecessarily large input signals. In contrast, the new controller~\eqref{eq:ICFC} is able to satisfy the input constraints and achieves a comparable performance of the tracking error.

\section{Conclusion}\label{Sec:Concl}

In the present paper we proposed a new funnel controller for a large class of nonlinear systems modelled by functional differential equations in the presence of input constraints. The funnel control law incorporates a novel dynamic adaptation scheme for the funnel boundaries, the asymptotic shape of which can be prescribed by the choice of controller design parameters, but which are widened according to the dynamics whenever the input saturation is active. We have rigorously proved that this controller achieves the control objective and, for the system class from the recent work~\cite{BergIlch21}, all involved signals are bounded for sufficiently large saturation level. The new controller extends the earlier funnel controller design from~\cite{BergIlch21} and shows a favorable performance in the presence of input constraints in comparative simulations.

Simulations have also shown that the controller performance strongly depends on the choice of the design parameters~\eqref{eq:FC-param}. In particular, the surjection~$N$ seems to influence the results a lot. Further research should reveal which choices should be preferred over others.




\appendices

\section{Proof of Theorem~\ref{Thm:FunCon}}\label{app:proof_main}

The proof consists of several steps.

\emph{Step 1}: We recast the closed-loop system in the form of an initial-value problem to which a well-known existence theory applies. First define
\begin{align*}
  \cD_1 &:= \setdef{(\eta,\psi)\in\R^m\times\R}{\|\eta\| < \psi},\\
  \gamma &: \cD_1\to\R^m,\ (\eta,\psi) \mapsto \left(1-\frac{\|\eta\|^2}{\psi^2}\right)^{-1} \eta.
\end{align*}
Next, we introduce continuous maps $\rho_k\colon\cD_k\to \R^m$, $k=1,\ldots,r$, recursively as follows:
\begin{equation*}\label{eq:rho-k}
\begin{aligned}
&\rho_1:\cD_1\to\R^m,\ (\eta_1,\psi_1)\mapsto \eta_1,\\
&\cD_k:= \Bigg\{(\eta_1,\ldots,\eta_k,\psi_1,\ldots,\psi_k)\in\R^{km+k}\ \Bigg| \\
&\!\!\!\!\begin{array}{l} (\eta_1,\ldots,\eta_{k-1},\psi_1,\ldots,\psi_{k-1})\in\cD_{k-1},\\
 \big(\eta_k\!+\!\gamma(\rho_{k-1}(\eta_1,\ldots,\eta_{k-1},\psi_1,\ldots,\psi_{k-1}),\psi_{k-1}),\psi_k\big) \!\in\!\cD_1\end{array}\!\!\!\!\Bigg\}\!\,,\\
&\rho_k\colon\cD_k\to\R^m,(\eta_1,\ldots,\eta_k,\psi_1,\ldots,\psi_k)\mapsto \eta_k \\
&\qquad\qquad+\gamma (\rho_{k-1}{ (\eta_1,\ldots,\eta_{k-1},\psi_1,\ldots,\psi_{k-1}),\psi_{k-1})}).
\end{aligned}
\end{equation*}
Note that each of the sets~$\cD_k$ is non-empty and open. Set~$n=r(m+1)$ and write $\xi\in\R^n$ as $(\xi_1,\ldots,\xi_{2r})$ with $\xi_i\in\R^m$ for $i=1,\ldots,r$ and $\xi_i\in\R$ for $i=r+1,\ldots,2r$. Define, for $k=1,\ldots,r$,
\begin{align*}
    \pi_k &: \R_{\ge 0}\times\R^n\to\R^{km+k},\\
     &(t,\xi)\mapsto \big(\xi_1-y_{\text{\rm ref}}(t),\ldots, \xi_k-y_{\text{\rm ref}}^{(k-1)}(t), \xi_{r+1},\ldots, \xi_{r+k}\big),\\
\cD&:=\setdef{(t,\xi)\in\R_{\ge 0}\times\R^n }{ \!\! \begin{array}{l} \pi_r (t,\xi)\in\cD_r,\\
 \xi_{r+i} > \frac{\beta_i}{\alpha_i} \text{ for $i=1,\ldots,r$}\end{array}\!\!\! },
\end{align*}
where $\cD$ is non-empty and relatively open. Further define
\begin{align*}
    V:\cD\to\R^m,\ (t,\xi)\mapsto N\left(\tfrac{1}{1-\tfrac{\|\rho_r(\pi_r(t,\xi))\|^2}{\xi_{2r}^2}}\right) \rho_r(\pi_r(t,\xi)).
\end{align*}
and
\begin{align*}
  &F\colon \cD\times\R^q\to\R^n,\\
  &(t,\xi,\eta)\mapsto \begin{pmatrix} \xi_2\\\vdots\\\xi_r\\f\big(d(t),\eta,\sat(V(t,\xi))\big)\\
   p_1 \xi_{r+2} - \alpha_1\xi_{r+1} + \beta_1 - p_1 \frac{\beta_{2}}{\alpha_{2}}\\ \vdots \\
   p_{r-1} \xi_{2r} - \alpha_{r-1}\xi_{2r-1} + \beta_{r-1} - p_{r-1} \frac{\beta_{r}}{\alpha_{r}}\\
   -\alpha_r \xi_{2r} +\beta_r + \xi_{2r} \frac{\|V(t,\xi)-\sat(V(t,\xi))\|}{\|\rho_r(\pi_r(t,\xi))\|}
   \end{pmatrix}.
\end{align*}
Note that the function~$F$, and in particular its last component, is well-defined on $\cD\times\R^q$: Since $N$ is continuous and $\xi_{2r}>\tfrac{\beta_r}{\alpha_r}$, there exists $\delta>0$ such that for all $(t,\xi)\in\cD$ with $\|\rho_r(\pi_r(t,\xi))\|< \delta$ we have that $\|V(t,\xi)\| < \theta$ for $\theta$ as in~(P5), and hence $\|V(t,\xi)-\sat(V(t,\xi))\| = 0$.

Writing
\[
x(t)=\big(y(t)^\top,\ldots,y^{(r-1)}(t)^\top, \psi_1(t),\ldots, \psi_r(t)\big)
\]
we see that the closed-loop initial-value problem~\eqref{eq:Sys},~\eqref{eq:IC},~\eqref{eq:ICFC}
may now be formulated as
\begin{equation}\label{eq:IVP-CL}
\begin{aligned}
\dot x(t)&=F\big(t,x(t),T(x)(t)\big),\\
x|_{[-h,0]}&=x^0\in C([-h,0],\R^n),
\end{aligned}
\end{equation}
where, for $t\in [-h,0]$,
\[
x^0(t):=\big( y^0(t)^\top,\ldots,(y^0)^{(r-1)}(t)^\top, \psi_1^0,\ldots, \psi_r^0\big)^\top.
\]
The function $F$ is measurable in~$t$, continuous in~$(\xi,\eta)$ and locally essentially bounded. By~\eqref{eq:InitCond} we see that $(0,x^0(0))\in\cD$. Therefore, an application of a variant of~\cite[Thm.~B.1]{IlchRyan09}\footnote{Although the property~(P3) of the operator~$T$ is weaker than required in~\cite{IlchRyan09}, this ``local'' property suffices for the proof.} yields the existence of a solution of~\eqref{eq:IVP-CL} and every solution can be extended to a maximal solution. Furthermore, any maximal solution $x:[-h,\omega)\to\R^{n}$, $\omega\in(0,\infty]$, of~\eqref{eq:IVP-CL} has the property that its graph
\[
    \cG := \setdef{(t,x(t))}{t\in[0,\omega)} \subset \cD
\]
has a closure which is not a compact subset of~$\cD$.

\emph{Step 2}: In this step we record some observations for later use. First observe that $(t,x(t))\in\cD$ for $t\in [0,\omega)$ implies that $\pi_i(t,x(t))\in\cD_i$ and hence we may define
\begin{align*}
    e_i(t) &:= \rho_i(\pi_i(t,x(t))),\quad \psi_i(t):= x_{r+i}(t),\\
    k_i(t)&:= \left(1- \frac{\|e_i(t)\|^2}{\psi_i(t)^2}\right)^{-1},\quad \gamma_i(t):= k_i(t) e_i(t),
\end{align*}
for $i=1,\ldots,r$, with which we arrive at the quantities in the control law~\eqref{eq:ICFC}; in particular $V(t,x(t)) = N(k_r(t)) e_r(t) = v(t)$. Clearly, we have $\|e_i(t)\|<\psi_i(t)$ for all $t\in[0,\omega)$ and $i=1,\ldots,r$. Finally, from~\eqref{eq:ICFC} it follows that for almost all $t\in[0,\omega)$ we have
\begin{equation}\label{eq:ODE-ei}
\begin{aligned}
  \dot e_i(t)&= e_{i+1}(t) - \gamma_i(t) + \dot \gamma_{i-1}(t),\quad i=1,\ldots,r-1,\\
  \dot e_r(t) &= e^{(r)}(t) + \dot \gamma_{r-1}(t),
\end{aligned}
\end{equation}
where $\gamma_0(t) := 0$. For brevity, set $\kappa_i := \beta_i - p_i \tfrac{\beta_{i+1}}{\alpha_{i+1}}$ for $i=1,\ldots,r-1$.

\emph{Step 3}: We show that $\psi_i(t) \ge \mu_i(0) e^{-\alpha_i t} + \tfrac{\beta_i}{\alpha_i}$ for all $t\in[0,\omega)$ and $i=1,\ldots,r$, where $\mu_i(\cdot)$ is defined in statement~(iii). By $(t,x(t))\in\cD$ we have that $\psi_i(t)>0$, thus $\dot \psi_r(t) \ge -\alpha_r \psi_r(t) + \beta_r$ and hence $\psi_r(t) \ge \mu_r(0) e^{-\alpha_r t} + \tfrac{\beta_r}{\alpha_r}$ for all $t\in [0,\omega)$. Then, inductively for $i=r-1,\ldots,1$,
\begin{align*}
    \dot \psi_i(t) &\ge p_i \mu_{i+1} (0) e^{-\alpha_{i+1} t} + p_i \tfrac{\beta_{i+1}}{\alpha_{i+1}} - \alpha_i \psi_i(t) +\beta_i - p_i \tfrac{\beta_{i+1}}{\alpha_{i+1}} \\
    &\ge - \alpha_i \psi_i(t) +\beta_i,
\end{align*}
from which the claim follows.

\emph{Step 4}: We show that $\tfrac{\psi_i}{\psi_r}\in L^\infty([0,\omega),\R)$ for $i=1,\ldots,r-1$. Define
\[
    M_{r-1} := \max \left\{ \frac{\psi_{r-1}^0}{\psi_{r}^0},\frac{p_{r-1}\beta_r + \beta_{r-1}\alpha_r}{\beta_r (\alpha_{i}-\alpha_r)} \right\}
\]
and, recursively for $i=r-2,\ldots,1$,
\[
    M_i :=  \max \left\{ \frac{\psi_{i}^0}{\psi_{r}^0}, \frac{1}{\alpha_i-\alpha_r}\left(p_i M_{i+1} + \frac{\beta_i \alpha_r}{\beta_r}\right) \right\}.
\]
We prove the statement by showing
\begin{equation}\label{eq:Psi-i-r-bd}
    \forall\, i=1,\ldots,r-1 \ \forall\, t\in[0,\omega):\  \frac{\psi_i(t)}{\psi_r(t)} \le M_i
\end{equation}
by induction over $i=r-1,\ldots,1$. For $i=r-1$, seeking a contradiction, assume there exists $t_1\in[0,\omega)$ such that $\tfrac{\psi_{r-1}(t_1)}{\psi_r(t_1)}>M_{r-1}$, then $t_1>0$ and
\[
    t_0:= \max\setdef{t\in[0,t_1)}{\tfrac{\psi_{r-1}(t)}{\psi_r(t)}=M_{r-1}}
\]
 is well-defined. Therefore, $\tfrac{\psi_{r-1}(t)}{\psi_r(t)}\ge M_{r-1}$ for all $t\in[t_0,t_1]$ and hence
\begin{align*}
  \ddt \frac{\psi_{r-1}(t)}{\psi_r(t)} &=\frac{\dot \psi_{r-1}(t) \psi_r(t) - \psi_{r-1}(t) \dot \psi_r(t)}{\psi_r(t)^2} \\
  &\stackrel{\eqref{eq:ICFC}}{=} p_{r-1} - \alpha_{r-1} \frac{\psi_{r-1}(t)}{\psi_r(t)} + \frac{\kappa_{r-1}}{\psi_r(t)} + \alpha_{r} \frac{\psi_{r-1}(t)}{\psi_r(t)}\\
   &\quad - \frac{\beta_{r}}{\psi_r(t)} \cdot \frac{\psi_{r-1}(t)}{\psi_r(t)} - \frac{\kappa(v(t))}{\|e_r(t)\|} \cdot \frac{\psi_{r-1}(t)}{\psi_r(t)}\\
  &\stackrel{\rm Step\,3}{\le} p_{r-1} - \underset{>0\ \text{by~\eqref{eq:FC-param}}}{ \underbrace{(\alpha_{r-1}-\alpha_r)}} \frac{\psi_{r-1}(t)}{\psi_r(t)} + \frac{\beta_{r-1} \alpha_r}{\beta_r}\\
  &\le p_{r-1} - (\alpha_{r-1}-\alpha_r) M_{r-1} + \frac{\beta_{r-1} \alpha_r}{\beta_r} \le 0
\end{align*}
for almost all $t\in [t_0,t_1]$, from which we infer
\[
    M_{r-1} < \frac{\psi_{r-1}(t_1)}{\psi_r(t_1)} \le \frac{\psi_{r-1}(t_0)}{\psi_r(t_0)} = M_{r-1},
\]
a contradiction. Now let $i\in\{1,\ldots,r-2\}$ and assume the statement~\eqref{eq:Psi-i-r-bd} is true for $i+1$. Again, assume there exists $t_1\in[0,\omega)$ such that $\tfrac{\psi_{i}(t_1)}{\psi_r(t_1)}>M_{i}$ and define
\[
    t_0:= \max\setdef{t\in[0,t_1)}{\tfrac{\psi_{i}(t)}{\psi_r(t)}=M_{i}}
\]
so that $\tfrac{\psi_{i}(t)}{\psi_r(t)}\ge M_{i}$ for all $t\in[t_0,t_1]$. Then
\begin{align*}
  \ddt \frac{\psi_{i}(t)}{\psi_r(t)}&= p_i \frac{\psi_{i+1}(t)}{\psi_r(t)} - \alpha_{i} \frac{\psi_{i}(t)}{\psi_r(t)} + \frac{\kappa_{i}}{\psi_r(t)} + \alpha_{r} \frac{\psi_{i}(t)}{\psi_r(t)}\\
   &\quad - \frac{\beta_{r}}{\psi_r(t)} \cdot \frac{\psi_{i}(t)}{\psi_r(t)} - \frac{\kappa(v(t))}{\|e_r(t)\|} \cdot \frac{\psi_{i}(t)}{\psi_r(t)}\\
  &\le p_i M_{i+1} - (\alpha_{i}-\alpha_r) M_i + \frac{\beta_i \alpha_r}{\beta_r} \le 0,
\end{align*}
which again yields a contradiction and completes the proof of~\eqref{eq:Psi-i-r-bd}.

\emph{Step 5}: We show that $\tfrac{\psi_i}{\psi_{i+1}}\in L^\infty([0,\omega),\R)$ for $i=1,\ldots,r-2$. Fix $i\in\{1,\ldots,r-2\}$, define
\[
    M := \max \left\{ \frac{\psi_{i}^0}{\psi_{i+1}^0},\frac{p_i \beta_{i+1} + \beta_{i}\alpha_{i+1}}{\beta_{i+1} (\alpha_{i}-\alpha_{i+1})} \right\}
\]
and assume there exists $t_1\in[0,\omega)$ such that $\tfrac{\psi_{i}(t_1)}{\psi_{i+1}(t_1)}>M$. Further define
\[
    t_0:= \max\setdef{t\in[0,t_1)}{\tfrac{\psi_{i}(t)}{\psi_{i+1}(t)}=M},
\]
so that $\tfrac{\psi_{i}(t)}{\psi_{i+1}(t)}\ge M$ for all $t\in[t_0,t_1]$. Then
\begin{align*}
  \ddt \frac{\psi_{i}(t)}{\psi_{i+1}(t)}&\stackrel{\eqref{eq:ICFC}}{=} p_i - \alpha_{i} \frac{\psi_{i}(t)}{\psi_{i+1}(t)} + \frac{\kappa_{i}}{\psi_{i+1}(t)} + \alpha_{i+1} \frac{\psi_{i}(t)}{\psi_{i+1}(t)}\\
   &\quad - p_{i+1} \frac{\psi_i(t) \psi_{i+2}(t)}{\psi_{i+1}(t)^2} - \frac{\kappa_{i+1} \psi_{i}(t)}{\psi_{i+1}(t)^2}\\
  &\stackrel{\rm Step\,3}{\le}  p_i - \underset{>0\ \text{by~\eqref{eq:FC-param}}}{ \underbrace{(\alpha_{i}-\alpha_{i+1})}} \frac{\psi_{i}(t)}{\psi_{i+1}(t)} + \beta_i \frac{\alpha_{i+1}}{\beta_{i+1}} \\
  &\quad - \frac{\psi_i}{\psi_{i+1}} \left(\frac{\beta_{i+1}}{\psi_{i+1}} + \frac{p_{i+1}}{\psi_{i+1}} \left(\psi_{i+2} - \frac{\beta_{i+2}}{\alpha_{i+2}}\right)\right)  \\
  &\stackrel{\rm Step\,3}{\le} p_i - (\alpha_{i}-\alpha_{i+1})M + \frac{\beta_i \alpha_{i+1}}{\beta_{i+1}} \le  0
\end{align*}
for almost all $t\in [t_0,t_1]$, from which we infer
\[
    M < \frac{\psi_{i}(t_1)}{\psi_{i+1}(t_1)} \le \frac{\psi_{i}(t_0)}{\psi_{i+1}(t_0)} = M,
\]
a contradiction.

\emph{Step 6}: We show
\begin{multline*}
    \forall\, i=2,\ldots,r:\ \Big( k_1,\ldots,k_{i-1}\in L^\infty([0,\omega),\R)\\ \implies \frac{\dot \gamma_{i-1}}{\psi_i} \in L^\infty([0,\omega),\R^m)\Big)
\end{multline*}
by induction over~$i$. For $i=2$ we have that
\begin{align*}
  &\dot \gamma_1(t) = 2 k_1(t)^2 \left(\frac{\|e_1(t)\|^2}{\psi_1(t)^3} \dot \psi_1(t) + \frac{e_1(t)^\top \dot e_1(t)}{\psi_1(t)^2}\right) e_1(t)\\
  &\quad + k_1(t)\dot e_1(t)\\
  &\stackrel{\eqref{eq:ICFC},\eqref{eq:ODE-ei}}{=} 2 k_1(t)^2 \bigg(-\frac{\|e_1(t)\|^2}{\psi_1(t)^3} \big(p_1 \psi_2(t)-\alpha_1\psi_1(t)+\kappa_1\big) e_1(t) \\
  & + \frac{\|e_1(t)\|^2}{\psi_1(t)^2}\big(e_2(t)\!-\!k_1(t) e_1(t)\big)\bigg) + k_1(t)\big( e_2(t) \!-\! k_1(t) e_1(t)\big)
\end{align*}
for almost all $t\in[0,\omega)$. Since $\tfrac{\|e_i(t)\|}{\psi_i(t)}<1$ for $i=1,2$ and $k_1(\cdot)$ is bounded by assumption, we find that
\begin{align*}
  \frac{\|\dot \gamma_1(t)\|}{\psi_2(t)} &\le 2 \|k_1\|_\infty^2\bigg( 2p_1 \!+\! \alpha_1 \frac{\psi_1(t)}{\psi_2(t)} \!+\! \frac{\kappa_1}{\psi_2(t)}  \!+\! \|k_1\|_\infty\frac{\psi_1(t)}{\psi_2(t)}\bigg)\\
   &\quad+ \|k_1\|_\infty \bigg(1 + \|k_1\|_\infty \frac{\psi_1(t)}{\psi_2(t)}\bigg)
\end{align*}
for almost all $t\in[0,\omega)$ and by Step~3 and Step~5 it follows that $\tfrac{\dot \gamma_{1}}{\psi_2}$ is bounded. Now assume the assertion is true for $i\in\{2,\ldots,r-1\}$. Then we have that
\begin{align*}
  &\dot \gamma_i(t) = 2 k_i(t)^2 \left(\frac{\|e_i(t)\|^2}{\psi_i(t)^3} \dot \psi_i(t) + \frac{e_i(t)^\top \dot e_i(t)}{\psi_i(t)^2}\right) e_i(t)\\
  &\quad + k_i(t) \dot e_i(t)\\
  &\stackrel{\eqref{eq:ICFC},\eqref{eq:ODE-ei}}{=} 2 k_i(t)^2 \bigg(\!\!-\frac{\|e_i(t)\|^2}{\psi_i(t)^3} \big(p_i \psi_{i+1}(t)\!-\!\alpha_i\psi_i(t)\!+\!\kappa_i\big) e_i(t) \\
  &\quad + \frac{\|e_i(t)\|^2}{\psi_i(t)^2}\big(e_{i+1}(t)-k_i(t) e_i(t)+ \dot \gamma_{i-1}(t)\big)\bigg) \\
  &\quad + k_i(t)\big(e_{i+1}(t)-k_i(t) e_i(t)+ \dot \gamma_{i-1}(t)\big)
\end{align*}
for almost all $t\in[0,\omega)$, by which
\begin{align*}
  &\frac{\|\dot \gamma_i(t)\|}{\psi_{i+1}(t)} \le 2 \|k_i\|_\infty^2\bigg( 2p_i  \!+\! \alpha_i \frac{\psi_i(t)}{\psi_{i+1}(t)}  \!+\! \frac{\kappa_i}{\psi_{i+1}(t)}  \!+\! \frac{\|\dot \gamma_{i-1}(t)\|}{\psi_{i+1}(t)}\\
  & + \! \|k_i\|_\infty\frac{\psi_i(t)}{\psi_{i+1}(t)}\bigg) \!+\! \|k_i\|_\infty \bigg(1 \!+\! \|k_i\|_\infty \frac{\psi_i(t)}{\psi_{i+1}(t)} \!+\! \frac{\|\dot \gamma_{i-1}(t)\|}{\psi_{i+1}(t)}\bigg)
\end{align*}
for almost all $t\in[0,\omega)$. By presupposition and Step~5 we have that
\[
    \frac{\dot \gamma_{i-1}}{\psi_{i+1}} = \frac{\dot \gamma_{i-1}}{\psi_{i}} \cdot \frac{\psi_i}{\psi_{i+1}} \in L^\infty([0,\omega),\R^m),
\]
and hence boundedness of $\tfrac{\dot \gamma_{i}}{\psi_{i+1}}$ follows from Steps~3 and~5.

\emph{Step 7}: We show that, if $\omega<\infty$, then for all $i=1,\ldots,r$ the continuous functions
\[
    \rho_i:[0,\omega)\to\R,\ t\mapsto \frac{\|\psi_i|_{[0,t]}\|_\infty}{\psi_i(t)}
\]
are bounded. Fix $i\in\{1,\ldots,r\}$ and $\tau\in[0,\omega)$. Seeking a contradiction, assume  there exists a strictly increasing sequence $(t_k)\to\omega$ such that $\rho_i(t_k)\to\infty$ for $k\to\infty$. Define, for $k\in\N$,
\begin{align*}
    s_k&:= \sup \setdef{t\in[0,t_k]}{ \psi_i(t) = \|\psi_i|_{[0,t_k]}\|_\infty}.
\end{align*}
Since $\psi_j(t) > \tfrac{\beta_j}{\alpha_j}$ for all $j=1,\ldots,r$ by Step~3 we find that $\dot \psi_i(t) \ge -\alpha_i \psi_i(t) + \beta_i$ for all $t\in[0,\omega)$ and hence
\begin{align*}
    \psi_i(t_k)& \ge e^{-\alpha_i (t_k-s_k)} \psi_i(s_k) + \frac{\beta_i}{\alpha_i} \big(1-e^{-\alpha_i (t_k-s_k)}\big) \\
    & \ge e^{-\alpha_i \omega} \psi_i(s_k),
\end{align*}
by which
\begin{align*}
  \rho_i(t_k) &= \frac{\psi_i(s_k)}{\psi_i(t_k)} \le e^{\alpha_i\omega} < \infty,
\end{align*}
a contradiction.

\emph{Step 8}: We show that, if $k_1,\ldots,k_{r-1}\in L^\infty([0,\omega),\R)$, then there exists $c\in L^\infty_{\loc}([0,\omega),\R)$ such that for all $t\in [0,\omega)$ we have
\[
    \|f\big(d(t),T(y,\dot y, \ldots,y^{(r-1)})(t),\sat(v(t))\big)\| \le c(t) \psi_r(t),
\]
and if $\omega<\infty$, then $c\in L^\infty([0,\omega),\R)$. By the sector bound property~(P4) we have that
\begin{align*}
   & \|f\big(d(t),T(y,\dot y, \ldots,y^{(r-1)})(t),\sat(v(t))\big)\| \\
    &\quad \le M_{1}\big(t,d(t),\sat(v(t))\big) \\
    &\qquad + \sum_{i=1}^{r} M_{i+1}\big(t,d(t),\sat(v(t))\big) \|y^{(i-1)}|_{[-h,t]}\|_\infty
\end{align*}
for all $t\in [0,\omega)$. Furthermore, with $K_i(t) := \|(y^0)^{(i-1)}|_{[-h,0]}\|_\infty + \|y_{\rm ref}^{(i-1)}|_{[0,t]}\|_\infty$, we have
\begin{align*}
  &\frac{\|y^{(i-1)}|_{[-h,t]}\|_\infty}{\psi_r(t)} \le \frac{1}{\psi_r(t)}\big( \|e^{(i-1)}|_{[0,t]}\|_\infty + K_i(t)\big)\\
  &\stackrel{\eqref{eq:ICFC}}{\le} \frac{1}{\psi_r(t)}\big( \|e_i|_{[0,t]}\|_\infty + \|k_{i-1} e_{i-1}|_{[0,t]}\|_\infty + K_i(t)\big)\\
  &\le \frac{1}{\psi_r(t)}\big( \|\psi_i|_{[0,t]}\|_\infty + \|k_{i-1}\|_\infty \|\psi_{i-1}|_{[0,t]}\|_\infty + K_i(t)\big) \\
  &\le  \frac{\psi_i(t)}{\psi_r(t)} \rho_i(t) + \|k_{i-1}\|_\infty \frac{\psi_{i-1}(t)}{\psi_r(t)} \rho_ {i-1}(t) + \frac{\alpha_r K_i(t)}{\beta_r}
\end{align*}
for all $t\in[0,\omega)$ and all $i=1,\ldots,r$, where $k_0:= 0$. Then, since $\tfrac{\psi_i}{\psi_r}$ and $\tfrac{\psi_{i-1}}{\psi_r}$ are bounded by Step~4 there exist $c_{i,1},c_{i,2}\ge 0$ such that
\[
    \frac{\|y^{(i-1)}|_{[-h,t]}\|_\infty}{\psi_r(t)} \le c_{i,1} \rho_{i-1}(t) + c_{i,2} \rho_i(t) + \frac{\alpha_r K_i(t)}{\beta_r} =: \tilde c(t).
\]
Define
\begin{align*}
    c(t)&:=  M_{1}\big(t,d(t),\sat(v(t))\big) \frac{\alpha_r}{\beta_r}\\
    &\quad + \sum_{i=1}^{r} M_{i+1}\big(t,d(t),\sat(v(t))\big) \tilde c(t),
\end{align*}
then, since $1\le \tfrac{\alpha_r}{\beta_r}\psi_r(t)$, we have that
\begin{align*}
    \|f\big(d(t),T(y,\dot y, \ldots,y^{(r-1)})(t),\sat(v(t))\big)\| \le c(t) \psi_r(t)
\end{align*}
and the claim follows from the observation that $c$ is locally essentially bounded since $\rho_1,\ldots,\rho_r,M_1,\ldots,M_{r+1},K_1,\ldots,K_r$ are continuous and $d$, $\sat(v)$ are bounded. Furthermore, if $\omega<\infty$ then~$c$ is bounded, because $\rho_1,\ldots,\rho_r$ are bounded by Step~7 and $M_1,\ldots,M_{r+1}$ and $K_1,\ldots,K_r$ are bounded by continuity (on $\R_{\ge 0}\times\R^p\times\R^m$ in case of the former).

\emph{Step 9}: We show that $k_i\in L^\infty([0,\omega),\R)$ for $i=1,\ldots,r-1$ or, equivalently, $\|e_i(t)\| \le \eps_i \psi_i(t)$ for all $t\in [0,\omega)$ and some $\eps_i\in (0,1)$, by induction over~$i$. Consider $i=1$ and choose $\eps_1\in(0,1)$ such that, invoking~\eqref{eq:InitCond},
\[
    \eps_1 > \max\left\{ \frac{\|e_1(0)\|}{\psi_1^0}, \frac{1}{p_1},
    \sqrt{1-\frac{1}{p_1^2 \frac{\alpha_1 \beta_2}{\beta_1 \alpha_2} + p_1 \alpha_1}} \right\}.
\]
Seeking a contradiction, assume there exists $t_1\in [0,\omega)$ such that $\|e_1(t_1)\| > \eps_1 \psi_1(t_1)$. Since $\|e_1(0)\| \le \eps_1 \psi_1(0)$,
\[
    t_0 := \max\setdef{t\in[0,t_1)}{\|e_1(t)\|= \eps_1 \psi_1(t)}
\]
is well-defined. Then, for all $t\in[t_0,t_1]$, we have
\begin{align*}
  \|e_1(t)\| \ge \eps_1 \psi_1(t)\quad\text{and}\quad k_1(t) \ge \frac{1}{1-\eps_1^2}
\end{align*}
and obtain that
\begin{align*}
  &\tfrac12 \ddt \|e_1(t)\|^2 \stackrel{\eqref{eq:ODE-ei}}{=} e_1(t)^\top \big( e_{2}(t) - k_1(t) e_1(t)\big)\\
  &\quad\le \big( -k_1(t) \|e_1(t)\| + \psi_2(t)\big) \|e_1(t)\| \\
  &\quad\stackrel{\eqref{eq:ICFC}}{=} \big(\eps_1 \dot \psi_1(t) - \eps_1 p_1 \psi_2(t) + \eps_1 \alpha_1 \psi_1(t) - \eps_1 \kappa_1 \\
  &\qquad - k_1(t) \|e_1(t)\| + \psi_2(t) \big) \|e_1(t)\| \\
  &\stackrel{\eps_1 p_1>1}{\le} \left(\eps_1 \dot \psi_1(t) + \alpha_1 \psi_1(t) + \eps_1 p_1 \frac{\beta_2}{\alpha_2} - \frac{\eps_1\psi_1(t)}{1-\eps_1^2}\right)  \|e_1(t)\| \\
  &\quad= \eps_1 \dot \psi_1(t)  \|e_1(t)\| -  \left( \xi_1 \psi_1(t) - \xi_2 \right)  \|e_1(t)\|,
\end{align*}
where $\xi_1 := -\alpha_1 +  \frac{\eps_1}{1-\eps_1^2} \ge -\alpha_1 +  \frac{1}{p_1(1-\eps_1^2)} > 0$ by choice of $\eps_1$ and $\xi_2:=\eps_1 p_1 \tfrac{\beta_2}{\alpha_2}\le p_1 \tfrac{\beta_2}{\alpha_2}$. By Step~3 we further have that
\begin{align*}
    \xi_1 \psi_1(t) - \xi_2 &\ge \left(-\alpha_1 +  \frac{1}{p_1(1-\eps_1^2)}\right) \frac{\beta_1}{\alpha_1} -  p_1 \frac{\beta_2}{\alpha_2} > 0\\
    &\iff\quad  \frac{1}{1-\eps_1^2} > p_1^2 \frac{\alpha_1 \beta_2}{\beta_1 \alpha_2} + p_1 \alpha_1
\end{align*}
which is satisfied by choice of $\eps_1$. Therefore,  $\tfrac12 \ddt \|e_1(t)\|^2 \le \eps_1 \dot \psi_1(t)  \|e_1(t)\|$ for almost all $t\in [t_0,t_1]$ and hence
\begin{align*}
  \|e_1(t_1)\| - \|e_1(t_0)\| &= \int_{t_0}^{t_1}  \tfrac12 \|e_1(t)\|^{-1} \ddt \|e_1(t)\|^2 {\rm d}t\\
  &\le \int_{t_0}^{t_1}  \eps_1 \dot\psi_1(t) {\rm d}t = \eps_1 \psi_1(t_1) - \eps_1 \psi_1(t_0),
\end{align*}
which yields the contradiction
\[
    0 = \eps_1 \psi_1(t_0) - \|e_1(t_0)\|  \le \eps_1 \psi_1(t_1) - \|e_1(t_1)\| < 0.
\]
Therefore, $k_1\in L^\infty([0,\omega),\R)$. Now assume that $k_1,\ldots,k_{i-1}\in L^\infty([0,\omega),\R)$ for some $i\in\{2,\ldots,r-1\}$. Let $c_i:= \left\|\tfrac{\dot\gamma_{i-1}}{\psi_i}\right\|_\infty$, which exists by Step~6, and choose $\eps_i\in(0,1)$ such that, invoking~\eqref{eq:InitCond},
\[
    \eps_i > \max\left\{ \frac{\|e_i(0)\|}{\psi_i^0}, \frac{1}{p_i}, \sqrt{1-\frac{1}{p_i^2 \frac{\alpha_i \beta_{i+1}}{\beta_i\alpha_{i+1}} + p_i (\alpha_i+c_i)}} \right\}.
\]
Similar to the above arguments, seeking a contradiction, assume there exists $t_1\in [0,\omega)$ such that $\|e_i(t_1)\| > \eps_i \psi_i(t_1)$ and define
\[
    t_0 := \max\setdef{t\in[0,t_1)}{\|e_i(t)\|= \eps_i \psi_i(t)}
\]
so that, for all $t\in[t_0,t_1]$, we have
\begin{align*}
  \|e_i(t)\| \ge \eps_i \psi_i(t)\quad \text{and}\quad
  k_i(t)  \ge  \frac{1}{1-\eps_i^2}.
\end{align*}
Then we obtain that
\begin{align*}
  &\tfrac12 \ddt \|e_i(t)\|^2 \stackrel{\eqref{eq:ODE-ei}}{=} e_i(t)^\top \big( e_{i+1}(t) - k_i(t) e_i(t) + \dot \gamma_{i-1}(t)\big)\\
  &\quad\le \big( -k_i(t) \|e_i(t)\| + \psi_{i+1}(t) + c_i \psi_i(t)\big) \|e_i(t)\| \\
  &\quad\stackrel{\eqref{eq:ICFC}}{=} \big(\eps_i \dot \psi_i(t) - \eps_i p_i \psi_{i+1}(t) + \eps_i \alpha_i \psi_i(t) - \eps_i \kappa_i\\
   &\qquad - k_i(t) \|e_i(t)\| + \psi_{i+1}(t) + c_i \psi_i(t) \big) \|e_i(t)\| \\
  &\stackrel{\eps_i p_i>1}{\le} \Big(\eps_i \dot \psi_i(t) + \alpha_i \psi_i(t) + \eps_i p_i \frac{\beta_{i+1}}{\alpha_{i+1}} - \frac{\eps_i\psi_i(t)}{1-\eps_i^2} \\
  &\qquad+ c_i \psi_i(t)\Big)  \|e_i(t)\|\\
  &\quad= \eps_i \dot \psi_i(t)  \|e_i(t)\| -  \left( \xi_{i,1} \psi_i(t) - \xi_{i,2} \right)  \|e_i(t)\|,
\end{align*}
where $\xi_{i,1} := -\alpha_i - c_i +  \frac{\eps_i}{1-\eps_i^2} \ge -\alpha_i - c_i +  \frac{1}{p_i(1-\eps_i^2)} > 0$ by choice of $\eps_i$ and $\xi_{i,2}:=\eps_i p_i \frac{\beta_{i+1}}{\alpha_{i+1}}\le p_i \frac{\beta_{i+1}}{\alpha_{i+1}}$. Furthermore, Step~3 gives
\begin{align*}
    \xi_{i,1} \psi_i(t) - \xi_{i,2} &\ge \left(-\alpha_i \!-\! c_i \!+\!  \frac{1}{p_i(1-\eps_i^2)}\right) \frac{\beta_i}{\alpha_i} \!-\!p_i \frac{\beta_{i+1}}{\alpha_{i+1}} > 0\\
    &\iff\quad  \frac{1}{1-\eps_i^2} > p_i^2 \frac{\alpha_i \beta_{i+1}}{\beta_i\alpha_{i+1}} + p_i (\alpha_i+c_i)
\end{align*}
which is satisfied by choice of $\eps_i$. Therefore,  $\tfrac12 \ddt \|e_i(t)\|^2 \le \dot \eps_i \psi_i(t)  \|e_i(t)\|$ for almost all $t\in [t_0,t_1]$ and as above for $i=1$ an integration yields
\[
    0 = \eps_i \psi_i(t_0) - \|e_i(t_0)\|  \le \eps_i \psi_i(t_1) - \|e_i(t_1)\| < 0,
\]
a contradiction.

\emph{Step 10}: We show that, if $\omega<\infty$, then $k_r\in L^\infty([0,\omega),\R)$. Let $R:=\sup_{t\in[0,\omega)} \|y_{\rm ref}^{(r)}(t)\|$ and $c_r:= \left\|\tfrac{\dot\gamma_{r-1}}{\psi_r}\right\|_\infty$, which exists by Steps~6 and~9, and observe that by Step~8 and $\omega<\infty$ we have that $c(\cdot)$ is bounded. Further note that, invoking~\eqref{eq:ICFC} we may estimate
\begin{equation}\label{eq:est-kappa}
   \forall\, t\in [0,\omega):\ \kappa(v(t)) \ge |N(k_r(t))| \cdot \|e_r(t)\| - M,
\end{equation}
where $M>0$ is some upper bound of $\sat$, i.e., $\|\sat(v)\|\le M$ for all $v\in\R^m$. Now choose
\begin{equation}\label{eq:delta}
    \delta > \alpha_r + \|c\|_\infty + c_r + \frac{\alpha_r}{\beta_r} (M+R)
\end{equation}
and $\eps_r\in (0,1)$ so that, invoking~\eqref{eq:InitCond},
\begin{align*}
    \eps_r > \frac{\|e_r(0)\|}{\psi_r^0}\quad\text{and}\quad \eps_r \left| N\left(\frac{1}{1-\eps_r^2}\right)\right| \ge 2\delta,
\end{align*}
where the latter is possible because of the properties of~$N$ in~\eqref{eq:FC-param}. We show that $\|e_r(t)\| \le \eps_r \psi_r(t)$ for all $t\in [0,\omega)$, which is equivalent to $k_r\in L^\infty([0,\omega),\R)$. Seeking a contradiction, assume there exists $t_1\in [0,\omega)$ such that $\|e_r(t_1)\| > \eps_r \psi_r(t_1)$ and define
\[
    t_0 := \max\setdef{t\in[0,t_1)}{\|e_r(t)\|= \eps_r \psi_r(t)}.
\]
Then, for all $t\in[t_0,t_1]$, we have
\begin{align}\label{eq:est-er-kr}
  \|e_r(t)\| \ge \eps_r \psi_r(t) 
  \quad \text{and}\quad
  k_r(t)  \ge  \frac{1}{1-\eps_r^2}.
\end{align}
Since $|N(k_r(t_0))| = |N(\tfrac{1}{1-\eps_r^2})| \ge 2\delta/\eps_r$, there exists $t_2\in (t_0,t_1]$ such that
\[
   \forall\, t\in [t_0,t_2]:\ |N(k_r(t))|\ge \frac{\delta}{\eps_r}.
\]
Furthermore, by definition of~$t_0$ we have that $\|e_r(t_2)\| > \eps_r \psi_r(t_2)$. Then we obtain that
\begin{align*}
  &\tfrac12 \ddt \|e_r(t)\|^2 \stackrel{\eqref{eq:ODE-ei}}{=} e_r(t)^\top \big( e^{(r)}(t) + \dot\gamma_{r-1}(t)\big)\\
  &\quad\stackrel{\eqref{eq:Sys}}{\le}  \big( \|f\big(d(t),T(y,\dot y, \ldots,y^{(r-1)})(t),\sat(v(t))\big)\| \\
  &\qquad + \|y_{\rm ref}^{(r)}(t)\| + \|\dot\gamma_{r-1}(t)\|\big) \|e_r(t)\|\\
  &\ \stackrel{\rm Step\,8}{\le}  \big( c(t) \psi_r(t) + R + c_r \psi_r(t)\big) \|e_r(t)\|\\
  &\quad\stackrel{\eqref{eq:ICFC}}{=} \Big(\eps_r \dot \psi_r(t) + \eps_r \alpha_r \psi_r(t) - \eps_r \beta_r - \eps_r \psi_r(t) \tfrac{\kappa(v(t))}{\|e_r(t)\|}\\
  &\qquad  + c(t) \psi_r(t) + R + c_r \psi_r(t) \Big) \|e_r(t)\|\\
  &\stackrel{\eqref{eq:est-kappa},\eqref{eq:est-er-kr}}{\le} \Big(\eps_r \dot \psi_r(t)  - \eps_r \beta_r + M + R- \big( \eps_r |N(k_r(t))| - \eps_r \alpha_r   \\
  &\qquad - \|c\|_\infty - c_r\big) \psi_r(t) \Big) \|e_r(t)\|\\
  &\quad \le \Big(\eps_r \dot \psi_r(t) \!+\! M \!+\! R \!-\! \underset{>0\ \text{by~\eqref{eq:delta}}}{\underbrace{\big( \delta \!-\! \alpha_r  \!-\! \|c\|_\infty \!-\! c_r\big)}} \psi_r(t) \Big) \|e_r(t)\|\\
  &\ \stackrel{\rm Step\,3}{\le} \Big(\eps_r \dot \psi_r(t) \!+\! M \!+\! R \!-\! \big( \delta \!-\! \alpha_r  \!-\! \|c\|_\infty \!-\! c_r\big) \tfrac{\beta_r}{\alpha_r} \Big) \|e_r(t)\|\\
  &\quad \stackrel{\eqref{eq:delta}}{\le} \eps_r \dot \psi_r(t) \|e_r(t)\|
\end{align*}
for almost all $t\in[t_0,t_2]$ and upon integration we get
\begin{align*}
  \|e_r(t_2)\| - \|e_r(t_0)\| &= \int_{t_0}^{t_2}  \tfrac12 \|e_r(t)\|^{-1} \ddt \|e_r(t)\|^2 {\rm d}t\\
  &\le \int_{t_0}^{t_2}  \eps_r \dot \psi_r(t) {\rm d}t = \eps_r \psi_r(t_2) - \eps_r \psi_r(t_0),
\end{align*}
which yields the contradiction
\[
    0 = \eps_r \psi_r(t_0) - \|e_r(t_0)\|  \le \eps_r \psi_r(t_2) - \|e_r(t_2)\| < 0.
\]

\emph{Step 11}: We show that $\omega=\infty$, i.e., assertion~(i) of the theorem. Suppose that $\omega<\infty$. From Steps~9 and~10 it follows that $k_1,\ldots,k_r\in L^\infty([0,\omega),\R)$ and hence there exists $\nu_1\in (0,1)$ such that
\[
    \forall\, i=1,\ldots,r\ \forall\, t\in[0,\omega):\ \|e_i(t)\|\le \nu_1 \psi_i(t).
\]
Furthermore, by Step~3 we have $\psi_i(t)\ge \mu_i(0) e^{-\alpha_i \omega} + \tfrac{\beta_i}{\alpha_i} > \tfrac{\beta_i}{\alpha_i}$ for all $t\in[0,\omega)$, and hence there exists $\nu_2 >0$ such that
\[
    \forall\, i=1,\ldots,r\ \forall\, t\in[0,\omega):\ \psi_i(t)\ge \tfrac{\beta_i}{\alpha_i} + \nu_2.
\]
Moreover, from boundedness of $k_r$ it follows that~$\kappa(v)$ is bounded and hence $\tfrac{\kappa(v)}{\|e_r\|}$ is bounded, since $\kappa(v)$ vanishes when $\|e_r\|$ is small enough, cf.\ Step~1. Therefore, it follows from~\eqref{eq:ICFC} that there exist some $d_{r,1}, d_{r,2}\ge 0$ such that $\psi_r(t) \le d_{r,1} e^{d_{r,2} t} \le d_{r,1} e^{d_{r,2} \omega}$ for all $t\in [0,\omega)$. Then, a successive solution of the differential equations for $\psi_i$ in~\eqref{eq:ICFC} yields similar bounds for them for $i=1,\ldots,r-1$. Hence there exists $\nu_3>0$ such that
\[
    \forall\, i=1,\ldots,r\ \forall\, t\in[0,\omega):\ \|e_i(t)\| < \psi_i(t) \le \nu_3.
\]
Define
\[
\hat \cD:= \setdef{(t,\xi)\in [0,\omega]\times\R^n}{\!\!\begin{array}{l} \|\rho_i(\pi_i(t,\xi))\|\le \nu_1 \xi_{r+i},\\
\tfrac{\beta_i}{\alpha_i} + \nu_2 \le \xi_{r+i} \le \nu_3\\
\text{for } i=1,\ldots,r\end{array}\!\!\!},
\]
which is evidently a compact subset of $\cD$ since $y_{\rm ref},\ldots,y_{\rm ref}^{(r-1)}$ are bounded on $[0,\omega]$. Since $(t,x(t))\in \hat \cD$ for all $t\in[0,\omega)$, it follows that the set~$\cG$ from Step~1 is a compact subset of~$\cD$, a contradiction. Therefore, $\omega=\infty$.

\emph{Step 12}: We complete the proof by establishing assertions~(ii) and~(iii) of the theorem. Assertion~(ii) is a consequence of Steps~1,~9 and~11. Let  $[t_0,t_1)\subseteq\R_{\ge 0}$ with $t_1\in (t_0,\infty]$ be an interval with $v(t) = \sat(v(t))$ for all $t\in [t_0,t_1)$, then we prove~(iii) by induction over $i=r,\ldots,1$. For $i=r$ the statement is clear since $\dot \psi_r(t) = -\alpha_r \psi_r(t) + \beta_r$ for all $t\in[t_0,t_1)$. Suppose the statement is true for some $i\in\{2,\ldots,r\}$, then we first observe that
\begin{align*}
    &\int_{t_0}^t e^{-\alpha_{i-1} (t-s)} e^{-\alpha_j (s-t_0)} {\rm d}s \\
    &= \tfrac{1}{\alpha_{i-1}-\alpha_j} \big( e^{-\alpha_j (t-t_0)} - e^{-\alpha_{i-1} (t-t_0)}\big)
    \le \tfrac{1}{\alpha_{i-1}-\alpha_j} e^{-\alpha_j (t-t_0)}
\end{align*}
for all $t\in[t_0,t_1)$ and all $j=i+1,\ldots,r$. Therefore, upon solving the differential equation for~$\psi_{i-1}$ in~\eqref{eq:ICFC} over $[t_0,t_1)$ we find that
\begin{align*}
    \psi_{i-1}(t) &= e^{-\alpha_{i-1}(t-t_0)} \psi_{i-1}(t_0) + \frac{\kappa_{i-1}}{\alpha_{i-1}} \big(1-e^{-\alpha_{i-1} (t-t_0)}\big)\\
    &\quad + \int_{t_0}^t p_{i-1} e^{-\alpha_{i-1} (t-s)}\psi_i(s) {\rm d}s \\
    &\le  \frac{\beta_{i-1}}{\alpha_{i-1}} + \mu_{i-1}(t_0) e^{-\alpha_{i-1}(t-t_0)} \\
    &\quad + \sum_{j=i}^r p_{i-1} \mu_j(t_0) \nu_{ij}  \int_{t_0}^t e^{-\alpha_{i-1} (t-s)} e^{-\alpha_j (s-t_0)} {\rm d}s \\
    &\le \frac{\beta_{i-1}}{\alpha_{i-1}} + \sum_{j=i-1}^r \mu_j(t_0) \nu_{i-1,j} e^{-\alpha_j (t-t_0)}
\end{align*}
for all $t\in[t_0,t_1)$. This completes the proof.\hfill $\Box$

\section{Proof of Theorem~\ref{Thm:FunCon-BIR}}\label{app:proof_Cor}

\emph{Step 1}: We provide a constructive definition of the constant~$M$. First observe that Steps~1--6 are the same as in the proof of Theorem~\ref{Thm:FunCon}, Steps~7 and~8 are not needed here and Step~9 is again the same as in the proof of Theorem~\ref{Thm:FunCon}. With this we arrive at a maximal  solution~$(y,\psi_1,\ldots,\psi_r):[-h,\omega)\to\R^{m+r}$, $\omega\in(0,\infty]$, of~\eqref{eq:Sys},~\eqref{eq:IC},~\eqref{eq:ICFC} with bounded $k_1,\ldots,k_{r-1}$. Next we seek to define explicit bounds for the latter. To this end, define
\[
    M_i := \max \left\{ \frac{\psi_{i}^0}{\psi_{i+1}^0},\frac{p_i \beta_{i+1} + \beta_{i}\alpha_{i+1}}{\beta_{i+1} (\alpha_{i}-\alpha_{i+1})} \right\}
\]
for $i=1,\ldots,r-1$ and observe that by Steps~4 and~5 in the proof of Theorem~\ref{Thm:FunCon} we have that
\[
    \forall\, t\in [0,\omega):\ \frac{\psi_i(t)}{\psi_{i+1}(t)} \le M_i.
\]
Then, recursively define
\begin{align*}
    c_1 &:= 0,\\
    \eps_1 &:= \max\left\{ \eps, \frac{1}{p_1},
    \sqrt{1-\frac{1}{p_1^2 \frac{\alpha_1 \beta_2}{\beta_1 \alpha_2} + p_1 \alpha_1}} \right\},
\end{align*}
and for $i=2,\ldots,r$
\begin{align*}
     c_i &:= \frac{2}{1-\eps_{i-1}^2}\left(2p_{i-1} + \alpha_{i-1} M_{i-1} + \frac{\alpha_{i-1}\kappa_{i-1}}{\beta_{i-1}} + M_{i-1} c_{i-1} \right.\\
     &\quad\left.+ \frac{M_{i-1}}{1-\eps_{i-1}^2}\right) + \frac{1}{1-\eps_{i-1}^2} \left(1+ \frac{M_{i-1}}{1-\eps_{i-1}^2} + M_{i-1} c_{i-1}\right),\\
     \eps_i &:= \max\left\{ \eps, \frac{1}{p_i}, \sqrt{1-\frac{1}{p_i^2 \frac{\alpha_i \beta_{i+1}}{\beta_i\alpha_{i+1}} + p_i (\alpha_i+c_i)}} \right\}.
\end{align*}
It follows from Steps~6 and~9 in the proof of Theorem~\ref{Thm:FunCon} that
\[
    \frac{\|\dot \gamma_{i-1}(t)\|}{\psi_i(t)} \le c_i\quad\text{and}\quad k_{i-1}(t)\le\frac{1}{1-\eps_{i-1}^2}
\]
for all $t\in[0,\omega)$ and all $i=2,\ldots,r$.

Now we will use the high-gain property of~$f$ from property~(P6). Choose $K_p := \setdef{\delta\in\R^p}{\|\delta\| \le \|d\|_\infty}$ and define
\[
    \psi_i^{\max} := \frac{\beta_i}{\alpha_i} + \sum_{j=i}^r \nu_{ij} \left(\psi_i^0 - \frac{\beta_i}{\alpha_i}\right).
\]
for $i=1,\ldots,r$, where $\nu_{ij}$ is given in statement~(iii) of Theorem~\ref{Thm:FunCon},
\[
    B\!:=\! \setdef{\!\zeta\!\in\! C([-h,\infty),\R^{rm})}{\!\!\!\begin{array}{l} \|\zeta_1\|_\infty \le \psi_1^{\max} \!+\! \hat K,\\ \|\zeta_i\|_\infty\le \psi_{i}^{\max} \!+\! \tfrac{\psi_{i-1}^{\max}}{1-\eps_{i-1}^2} \!+\! \hat K,\\ i=2,\ldots,r\end{array}\!\!\!}\!,
\]
where $\hat K:= K + \max_{i=0,\ldots,r-1} \|(y^0)^{(i)}\|_\infty$ (recall that by assumption $K>0$ is given such that $\|y_{\rm ref}^{(i)}\|_\infty \le K$ for $i=0,\ldots,r$), and
\[
    K_q := \setdef{z\in\R^q}{ \|z\|\le \sup_{\zeta\in B} \| T(\zeta)\|_\infty},
\]
which is a compact set since~$T$ satisfies~(P3) for $\tau = \infty$.
Further set $\nu^*:= \tfrac12$, for which we obtain the corresponding function~$\chi$ as in~(P6).

Next we define
\[
    \chi^* := \big( K + c_r\psi_{r}^{\max}  + \alpha_r \psi_{r}^{\max}\big) \psi_{r}^{\max}
\]
and choose $\eps_r\in [\eps,1)$ such that
\[
    \chi\left(N\left(\frac{1}{1-\eps_r^2}\right)\right) \ge 2 \chi^*,
\]
which is possible because of the properties of~$\chi$ and~$N$. Finally, we define
\[
    M:= \left( \sup_{s\in \left[0,\tfrac{1}{1-\eps_r^2}\right]} |N(s)|\right) \psi_r^0 + \delta
\]
for some arbitrary $\delta>0$. Now let the saturation function $\sat$ be such that it satisfies~(P5) with $\theta=M$ (and note that the above derived properties of the solution are indeed independent of this property of $\sat$).

\emph{Step 2}: We show that $\|v(t)\|< M$ for all $t\in[0,\omega)$. Define
\[
    t_0 :=\inf \setdef{t\in [0,\omega)}{ \|v(t)\| \ge M}.
\]
Seeking a contradiction, assume that $t_0 < \omega$. By assumption we have that $\frac{\|e_r(0)\|}{\psi_r(0)} \le \eps$, thus $k_r(0) \le \frac{1}{1-\eps^2}\le \frac{1}{1-\eps_r^2}$. Therefore, we find that
\[
    \|v(0)\| \le |N(k_r(0))| \psi_r(0) < M,
\]
and hence $t_0>0$ and we have that $\|v(t)\|\le M$ for all $t\in [0,t_0]$. This also implies that $\sat(v(t)) = v(t)$ and hence $\kappa(v(t))=0$ for all $t\in [0,t_0]$. Then it follows from statement~(iii) in Theorem~\ref{Thm:FunCon} that
\[
    \psi_i(t) \le \frac{\beta_i}{\alpha_i} + \sum_{j=i}^r \nu_{ij} \left(\psi_i^0 - \frac{\beta_i}{\alpha_i}\right) e^{-\alpha_j t} \le \psi_i^{\max}
\]
for all $t\in[0,t_0]$ and all $i=1,\ldots,r$. With this we find that by~\eqref{eq:ICFC}
\[
    \|y(t)\| \le \|e(t)\| + K \le \psi_{1}^{\max} + K
\]
and
\begin{align*}
   \|y^{(i)}(t)\| &\le  \|e^{(i)}(t)\| + K \le \|e_{i+1}(t)\|+ k_i(t) \|e_i(t)\| + K \\
   &\le \psi_{i+1}^{\max} + \tfrac{\psi_{i}^{\max}}{1-\eps_{i}^2} + K
\end{align*}
for all $t\in[0,t_0]$ and $i=1,\ldots,r-1$. Therefore, $\zeta\in C([-h,\infty),\R^{rm})$ defined by
\[
    \zeta_i(t) = \begin{cases} (y^0)^{(i-1)}(t), & t\in[-h,0],\\  y^{(i-1)}(t), & t\in[0,t_0],\\  y^{(i)}(t_0), & t\ge t_0\end{cases}
\]
satisfies $\zeta\in B$ and $T(y,\ldots,y^{(r-1)})(t) = T(\zeta)(t)$ for all $t\in [0,t_0]$, by which $T(y,\ldots,y^{(r-1)})(t)\in K_q$ for all $t\in [0,t_0]$.

\emph{Step~2a}: We show that $\|e_r(t)\| \le \eps_r \psi_r(t)$ for all $t\in [0,t_0]$. Seeking a contradiction, assume there exists $t_2\in [0,t_0]$ such that $\|e_r(t_2)\| > \eps_r \psi_r(t_2)$ and define
\[
    t_1 := \max\setdef{t\in[0,t_2)}{\|e_r(t)\|= \eps_r \psi_r(t)}.
\]
Then, for all $t\in[t_1,t_2]$, we have
\begin{align*}
  \|e_r(t)\| \ge \eps_r \psi_r(t) 
  \quad \text{and}\quad
  k_r(t)  \ge  \frac{1}{1-\eps_r^2}.
\end{align*}
Since $k_r(t_1) = \tfrac{1}{1-\eps_r^2}$ we find that
\[
    \chi\big(N(k_r(t_1))\big) \ge 2 \chi^*,
\]
hence there exists $t_3\in (t_1,t_2]$ such that
\[
   \forall\, t\in [t_1,t_3]:\ |\chi\big(N(k_r(t))\big)|\ge \chi^*.
\]
By definition of~$\chi$ we find that
\begin{multline*}
    e_r(t)^\top f\big(d(t),T(y,\dot y, \ldots,y^{(r-1)})(t), -k_r(t) e_r(t) \big) \\
    \le - \chi\big(N(k_r(t)\big) \le -\chi^*
\end{multline*}
for all $t\in [t_1,t_3]$ and since $\sat(v(t)) = v(t) = -k_r(t) e_r(t)$ it follows that
\begin{align*}
  &\tfrac12 \ddt \|e_r(t)\|^2 \stackrel{\eqref{eq:ODE-ei}}{=} e_r(t)^\top \big( e^{(r)}(t) + \dot\gamma_{r-1}(t)\big)\\
  &\quad\stackrel{\eqref{eq:Sys}}{\le}   e_r(t)^\top f\big(d(t),T(y,\dot y, \ldots,y^{(r-1)})(t), \sat(v(t)) \big) \\
  &\qquad  + \big(\|y_{\rm ref}^{(r)}(t)\| + \|\dot\gamma_{r-1}(t)\|\big) \|e_r(t)\|\\
  &\quad\le  -\chi^* + K \psi_{r}^{\max} + c_r \big(\psi_{r}^{\max}\big)^2\\
  &\quad = \eps_r \dot \psi_r(t) \|e_r(t)\| - \eps_r \big(-\alpha_r \psi_r(t) + \beta_r\big) \|e_r(t)\| \\
  &\qquad -\chi^* + K \psi_{r}^{\max} + c_r \big(\psi_{r}^{\max}\big)^2\\
  &\quad\stackrel{\eps_r < 1}{\le} \eps_r \dot \psi_r(t) \|e_r(t)\| + \alpha_r \big(\psi_{r}^{\max}\big)^2 \\
  &\qquad -\chi^* + K \psi_{r}^{\max} + c_r \big(\psi_{r}^{\max}\big)^2\\
  &\quad\le \eps_r \dot \psi_r(t) \|e_r(t)\|,
\end{align*}
where in the last step the definition of~$\chi^*$ was used. Then, using $\|e_r(t_3)\| > \eps_r \psi_r(t_3)$, similar to Step~10 in the proof of Theorem~\ref{Thm:FunCon}, a contradiction follows.

\emph{Step~2b}: We conclude the argument of Step~2. Since $k_r(t)  \le  \frac{1}{1-\eps_r^2}$ for all $t\in[0,t_0]$ by Step~2a, it follows that
\[
    \|v(t)\|  \le \left( \sup_{s\in \left[0,\tfrac{1}{1-\eps_r^2}\right]} |N(s)|\right) \psi_r^0 = M - \delta < M
\]
for all $t\in [0,t_0]$, which contradicts $\|v(t_0)\| \ge M$. Therefore, $\|v(t)\|< M$ for all $t\in[0,\omega)$.

\emph{Step 3}: Similar to Step~11 in the proof of Theorem~\ref{Thm:FunCon} it can be shown that $\omega = \infty$. Furthermore, the statements~(i) and~(ii) of Theorem~\ref{Thm:FunCon-BIR} have been shown above, where $\limsup_{t\to\infty} \psi_i(t) \le \frac{\beta_i}{\alpha_i}$ follows from the estimate in statement~(iii) in Theorem~\ref{Thm:FunCon}. This concludes the proof. \hfill $\Box$

\section{Proof of Corollary~\ref{Cor:FunCon-Inv}}\label{app:proof_Cor_2}

The proof is a straightforward modification of that of Theorem~\ref{Thm:FunCon-BIR} by replacing $\eps$ with~$0$ and taking the resulting expressions for~$\eps_i$ as the definition for~$\hat \eps_i$. The assertions then follow along the lines of the proofs of Theorem~\ref{Thm:FunCon-BIR} and Theorem~\ref{Thm:FunCon}, resp. \hfill $\Box$


\bibliographystyle{IEEEtran}

\begin{IEEEbiography}[{\includegraphics[width=1.05in,clip,keepaspectratio]{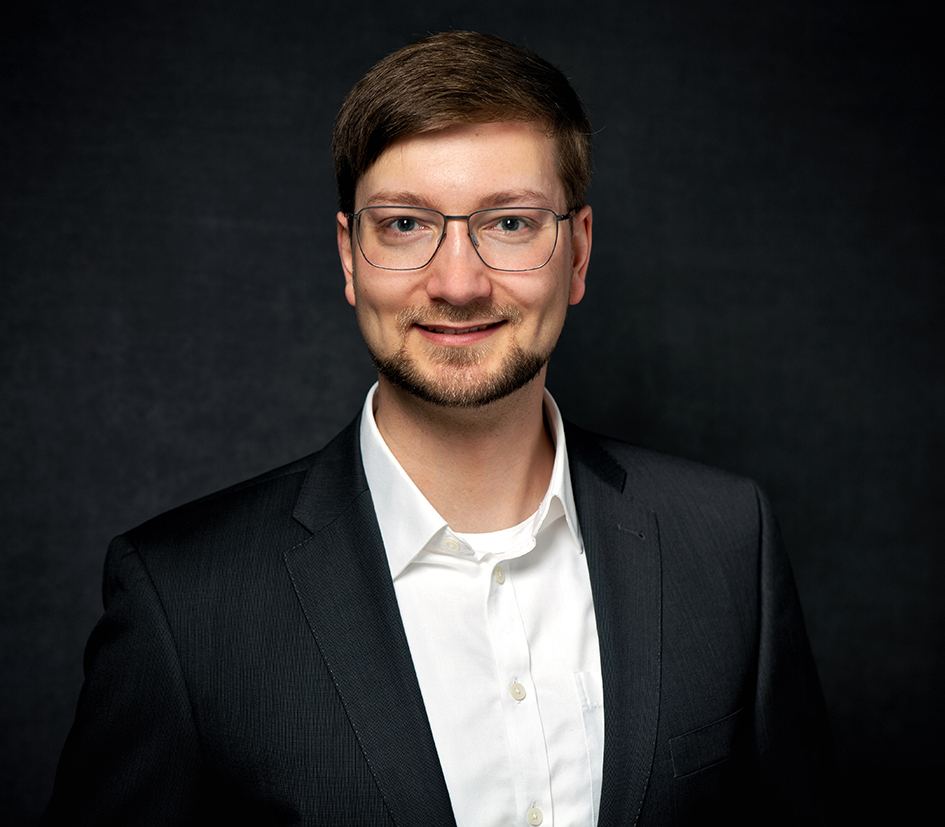}}]{Thomas Berger} was born in Germany in 1986. He received his B.Sc. (2008), M.Sc. (2010), and Ph.D. (2013), all in Mathematics and from Technische Universit\"at Ilmenau, Germany. From 2013 to 2018 Dr. Berger was a postdoctoral researcher at the Department of Mathematics, Universit\"at Hamburg, Germany. Since January 2019 he is a Juniorprofessor at the Institute for Mathematics, Universit\"at Paderborn, Germany. His research interest encompasses systems and control theory, differential–algebraic systems and multibody dynamics.

For his exceptional scientific achievements in the field of Applied Mathematics and Mechanics, Dr. Berger received the ``Richard-von-Mises Prize 2021'' of the International Association of Applied Mathematics and Mechanics (GAMM). He further received several awards for his dissertation, including the ``2015 European Ph.D. Award on Control for Complex and Heterogeneous Systems'' from the European Embedded Control Institute and the ``Dr.-Körper-Preis 2015'' from the GAMM. He serves as an Associate Editor for Mathematics of Control, Signals, and Systems, the IMA Journal of Mathematical Control and Information and the DAE Panel, and as a Review Editor for Frontiers in Control Engineering.
\end{IEEEbiography}

\end{document}